\documentclass[preprint,12pt]{elsarticle}
\usepackage{amssymb}
\usepackage{amsmath}
\usepackage{amsthm}
\usepackage{color}
\usepackage{float}
\usepackage[hidelinks]{hyperref}
\usepackage{booktabs}
\usepackage{bm}
\usepackage{xcolor}
\usepackage{tikz}
\usetikzlibrary{arrows.meta, positioning, calc, backgrounds, shadows}
\usepackage{xcolor}
\usepackage{amsmath, amssymb}
\usepackage{comment}
\usetikzlibrary{arrows.meta, positioning, calc, backgrounds, shadows}

\definecolor{myGreen}{rgb}{0.10, 0.67, 0.52}
\definecolor{myBlue}{rgb}{0.00, 0.45, 0.74}
\definecolor{myOrange}{rgb}{0.85, 0.55, 0.25}

\colorlet{bgGreen}{myGreen!7}
\colorlet{bgBlue}{myBlue!7}
\colorlet{bgOrange}{myOrange!10}

\journal{Physica D: Nonlinear Phenomena}

\begin{document}

\begin{frontmatter}

%\title{Stochastic separatrix and transition times in a temperature--phytoplankton model: a committor-based geometry}
\title{Geometric early warning indicator from stochastic separatrix structure in a random two-state ecosystem model}

\author[label1]{Yuzhu Shi}
\author[label1]{Larissa Serdukova}
\author[label2]{Yayun Zheng}
\author[label1,label3]{Sergei Petrovskii}
\author[label1,label4]{Valerio Lucarini}
%\author[label6]{J\"{u}rgen Kurths}

\affiliation[label1]{organisation={Department of Computing and Mathematical Sciences, University of Leicester,}, city={Leicester},
postcode={LE1 7RH},
country={UK}}
\affiliation[label2]{organisation={School of Mathematical Sciences, Jiangsu University,},
city={Zhenjiang},
postcode={212013},
country={China}}
\affiliation[label3]{organisation={Faculty of Physics, Mathematics and Natural Sciences, Peoples' Friendship University of Russia,},
city={Moscow},
postcode={117198},
           country={Russia}}
\affiliation[label4]{organisation={School of Sciences, Great Bay University,},
city={Dongguan},
postcode={523000},
country={China}}

\begin{abstract}
%\rev{
Under-ice blooms in the Arctic can develop rapidly under conditions where conventional early warning signals based on critical slowing down fail due to strong noise or limited observational records. We analyse noise-induced transitions in a temperature–phytoplankton stochastic differential equation model exhibiting bistability between background and bloom states. The committor function defines a stochastic separatrix as its $1/2$-isocommittor, and the normal width of the associated transition layer yields a geometric indicator via arc-length averaging. Under systematic variation of noise intensity, this indicator scales linearly with noise strength, while the logarithm of the mean first passage time follows the Freidlin–Wentzell asymptotic law. Eliminating the noise parameter produces an affine scaling between the logarithmic transition time and the inverse square of the geometric indicator. The relation is robust under variations in discretisation, neighbourhood definition, and diffusion structure, and holds in the weak-noise regime where the transition-layer width scales linearly with noise strength. Unlike variance or lag-one autocorrelation, the geometric indicator remains well defined when rapid transitions preclude reliable time-series estimation. These results provide a geometrically interpretable precursor of bloom onset that may support model-based ecological monitoring in high-variability Arctic systems.
%}
%Under-ice blooms in the Arctic can develop rapidly under conditions where conventional early warning signals based on critical slowing down fail due to strong noise or short observational records. We analyze noise induced transitions in a temperature--phytoplankton stochastic differential equation exhibiting bistability between a background biomass state and a bloom state. The committor function defines a stochastic separatrix as the isocommittor, and the normal width of the probabilistic transition layer along this curve yields a geometric indicator through arc-length averaging. Under systematic variation of noise intensity at fixed control parameter, the geometric indicator scales linearly with noise intensity while the logarithm of the mean first passage time follows the Freidlin--Wentzell asymptotic law. Eliminating the noise parameter couples these relations and produces an affine dependence of the logarithmic transition time on the inverse square of the geometric indicator. The scaling persists under variations in numerical discretisation, neighborhood size, and noise structure. Unlike variance or lag-one autocorrelation, the geometric indicator remains defined when rapid transitions prevent reliable estimation of time series statistics. This relation arises from the asymptotic structure of bistable diffusions with separable noise intensity and holds in the weak noise regime where the transition layer width scales linearly with the noise strength.
\end{abstract}

\end{frontmatter}

% =============================
\section{Introduction}
Since the late 1970s, Arctic sea ice has thinned continuously, and the melt season has lengthened markedly under global warming driven primarily by anthropogenic greenhouse gas emissions~\cite{NotzStroeve2016}. Observations show an earlier melt onset and a delayed autumn freeze-up, extending the melt season and enhancing summer solar radiation absorption through the ice-albedo feedback~\cite{Markus2009JGR,Stroeve2014MeltSeason}. As sea ice thins and the prevalence of melt ponds and leads increases, more solar radiation penetrates the ice to enter the upper ocean. Concurrently, surface layer freshening from ice melt strengthens stratification and suppresses vertical mixing, retaining phytoplankton within the euphotic zone for longer periods~\cite{Arrigo2012Science,ArdynaArrigo2020Frontiers}. In nutrient-replete shelf seas where biological losses from grazing and sinking remain low, phytoplankton biomass can accumulate rapidly over a few days to form under-ice blooms (UIBs)~\cite{Arrigo2014DSR2,Assmy2017SciRep}. The repeated observation of UIBs beneath continuous sea ice far from the ice edge, including during the 2011 Chukchi Sea cruise and the N--ICE2015 campaign, indicates that sea ice cover does not necessarily preclude high primary productivity. Under optimal light transmission and stratification, biomass can increase rapidly, substantially altering local biogeochemical cycling and food web structures~\cite{Arrigo2012Science,Lewis2020Science}.

However, rapidly proliferating UIBs are frequently accompanied by toxic harmful algal blooms (HABs), which place pressure on polar ecosystems and biological safety~\cite{Arrigo2012Science,ArdynaArrigo2020Frontiers,Arrigo2014DSR2,Anderson2012ARMS}. For example, monitoring in Arctic Alaska indicates that toxic phytoplankton expanding under warming and their toxins, including saxitoxins and domoic acid, exacerbate HAB risks. In the post-bloom period, toxins can persist and accumulate within food webs, adversely affecting higher trophic levels and fisheries resources~\cite{Anderson2022Oceanography,Lefebvre2016HarmfulAlgae,Gobler2020HarmfulAlgae}.

From a dynamical systems perspective, these abrupt events are characterised as regime shifts between two stable states, namely the background state and the bloom state~\cite{Scheffer2001,Lenton2008,PisarchikFeudel2014}.Transitions between these states are typically abrupt and exhibit hysteresis, making them difficult to reverse. Identifying early warning signals (EWS) prior to the development of HABs is therefore critical, as it provides essential lead time for toxin monitoring and emergency intervention~\cite{Anderson2012ARMS,Anderson2022Oceanography, Scheffer2009,Dakos2015}.

Over the past two decades, EWS have relied predominantly on time series statistical features derived from the critical slowing down (CSD) phenomenon. As a system approaches a bifurcation point under slow forcing, the decreasing local recovery rate manifests as increased variance and lag-one autocorrelation~\cite{Scheffer2009,Dakos2015}. However, their reliability is sensitive to sample length, detrending methods, and sampling frequency, and false alarms or missed detections are common under strong noise or in scenarios with random switching among different stable states~\cite{Boettiger2012,Laitinen2021ProbabilisticEWS,Morr2024}.

Additionally, the effectiveness of EWS is strongly affected by the choice of the observable used to construct them~\cite{Lohmannetal2025}. Recent theoretical developments show that a general theory of CSD in the vicinity of tipping points requires a multivariate rather than purely univariate statistical perspective, since tipping behaviour is often associated with changes in the geometry of the underlying attractor and its interaction with the basin boundary~\cite{SantosLucarini2022,LucariniChekroun2023,LucariniBodai2020}. From this viewpoint, indicators based on a single observable may fail to capture relevant dynamical information in high-dimensional systems. A broader discussion of EWS in ecological and climate contexts can be found in~\cite{Hastingsetal2026}. For Arctic UIBs, these challenges are compounded by the intrinsic episodic variability of light and stratification. Furthermore, sea ice cover obstructs satellite remote sensing, resulting in short, irregular observational records that hinder the reliable identification of indicators based on time series trends~\cite{Arrigo2012Science,ArdynaArrigo2020Frontiers}.

In recent years, several studies have explored early warning methods based on probabilistic structure and boundary geometry. In deterministic systems, bistability is governed by well-defined basin boundaries~\cite{Menck2013}. Under stochastic perturbations, noise reshapes the stochastic basins of attraction, replacing sharp deterministic boundaries with diffuse transition layers~\cite{Serdukova2016}. While initial conditions in a deterministic setting lead to unique attractors, noise renders these outcomes probabilistic. Analysing which phase-space regions approach the transition layer, and how this region shifts or widens with parametric changes, provides a robust metric for transition risk~\cite{JacquesDumas2023,JacquesDumas2024AIES,Kang2024NatCompSci}.

Motivated by these insights, we employ the committor function and the stochastic separatrix from transition path theory~\cite{WeinanE2006,Lucente2022}. For the UIB system, the committor function provides the conditional probability of reaching the bloom state before returning to the background state. Its isocommittor defines the stochastic separatrix~\cite{AntoniouSchwartz2009JCP}.Near this separatrix, the spatial gradient of the committor quantifies how rapidly transition probabilities change across state space. We thus propose a geometric early warning indicator defined by the transition layer width, and utilise the mean first passage time (MFPT) to quantify the time scale of transition difficulty~\cite{Arani2021ExitTime,Xu2023NonEqEWS}. In the weak-noise regime (Freidlin–Wentzell limit), the mean transition time is governed by the quasipotential barrier height, following a Kramers' type exponential law~\cite{Hanggi1990}.

We build upon the eco--climate coupled model proposed by Alsulami and Petrovskii~\cite{Alsulami2023}. By considering representative upper-ocean temperature and aggregate phytoplankton biomass, we construct a temperature--phytoplankton stochastic differential equation (SDE) by incorporating stochastic perturbations~\cite{Wouters2016,Franzke2015StochClimate,Berner2017StochParam}. Specifically, additive noise represents external thermal forcing fluctuations, while multiplicative noise is employed in the biomass equation to reflect disturbance intensity relative to biomass and ensure diffusion term positivity. To characterise bistability, we utilise a stochastic bifurcation diagram to examine the response of stationary statistics to parameters and noise intensity~\cite{Arnold1998Bifurcation,HomburgYoungGharaei2013,Risken1996}. We compute the committor function by solving a backward Kolmogorov equation and define the geometric indicator $EWS_{\mathrm{geom}}$ as the arc-length-averaged width of the transition layer around the stochastic separatrix $\Gamma = \{q = 1/2\}$. Separately, we compute the MFPT as a scalar measure of transition difficulty. 

%The main contributions of this work are threefold. First, we decompose the stochastic separatrix geometry into shift and transition layer width, converting the probability boundary into computable geometric indicators. Second, we analyze how these geometric quantities vary with parameters and noise intensity, and find that $EWS_{\mathrm{geom}}$ provides stable early warning signals and well-defined warning intervals under strong noise. Third, we derive an asymptotic scaling law linking geometry to time that in the weak-noise regime, the logarithmic MFPT linearly with the inverse square of $EWS_{\mathrm{geom}}$.

%\rev{
The novelty of this work lies primarily in the explicit geometric–temporal coupling between stochastic basin geometry and transition timescales. While quasipotential theory and transition path theory independently characterise rare-event transitions, we derive an explicit asymptotic relation linking the arc-length-averaged transition-layer width of the stochastic separatrix to the logarithmic MFPT. The resulting affine scaling $\log \langle \tau \rangle \sim 1/\mathrm{EWS}_{\mathrm{geom}}^2$ emerges from coupling two classical asymptotics through the common noise intensity parameter. This provides a geometric interpretation of transition difficulty in bistable diffusions with separable noise structure.
%}

The remainder of the paper is organised as follows. Section 2 starts from the deterministic model and introduces stochastic perturbations exhibiting state dependence to obtain the temperature--phytoplankton SDE. Section 3 defines the geometric decomposition of the stochastic separatrix into shift and transition-layer width $EWS_{\mathrm{geom}}$. Section 4 presents how these geometric quantities evolve with parameters and noise intensity, and compares its warning performance with classical time series indicators. Section 5 establishes the asymptotic relation between $EWS_{\mathrm{geom}}^2$ and MFPT in the weak-noise limit and analyses its validity conditions. Section 6 concludes and outlines future directions. For convenience, the main abbreviations used throughout the manuscript are collected in a glossary at the end of the paper.
% ===================================================================
\section{Temperature--phytoplankton model}\label{sec:model_methods}
\subsection{Deterministic model}\label{subsec:det_model}
To specifically investigate UIB dynamics, 
we consider a temperature–phytoplankton system consisting of a representative
upper-ocean environmental variable $T$, interpreted as an effective thermal/light
state under sea ice, and the aggregate phytoplankton biomass $u$ in the euphotic zone~\cite{Alsulami2023}:

\begin{align}
\frac{\mathrm{d}T}{\mathrm{d}t} &= \frac{1}{\gamma}\Big[-aT^{4}+b\big(1-S(u)\big)\Big], \label{eq:det_T}\\
\frac{\mathrm{d}u}{\mathrm{d}t} &= u\big(g(T)-\mu-u\big). \label{eq:det_u}
\end{align}

In Eq.~\eqref{eq:det_T}, the term $-aT^{4}$ accounts for Stefan--Boltzmann longwave radiative loss, where $\gamma$ is the effective thermal inertia. The term $b(1-S(u))$ represents the net shortwave absorption, where the function $S(u)$ parametrises the combined impact of sea ice cover, light transmission through the layer of ice, and water column absorption on the radiation budget. As biomass increases, pigment absorption reduces $S(u)$, thereby establishing a positive feedback loop. Phytoplankton dynamics follow a logistic growth law, where $g(T)$ is the temperature growth rate and $\mu$ represents the aggregate loss rate due to grazing, sinking, and viral lysis. The quadratic term $-u^2$ represents growth inhibition resulting from nutrient depletion and self-shading.

The functional forms of $S(u)$ and $g(T)$ are given by~\cite{Alsulami2023}:
\begin{equation}
S(u)=s_0+(s_1-s_0)\mathrm{e}^{-\alpha_1 u}, \qquad g(T)=b_1\exp\!\Big(-\frac{T_0}{T}-\alpha_2T\Big),
\label{eq:Su_gT}
\end{equation}
These formulations ensure $S'(u) < 0$, consistent with the physical expectation that increased biomass enhances heat absorption. The growth rate $g(T)$ is Arrhenius limited at low temperature and suppressed at high temperature. In this study, $b_1$ serves as the primary control parameter, interpreted as the effective growth rate that encapsulates variations in light, nutrients, and community composition. 
Model parameters are summarised in Table~\ref{tab:parameters}.

\begin{table}[H]
\centering \small
%\caption{Parameter values and descriptions.}
\caption{Numerical values of the dimensionless model parameters used in simulations~\cite{Alsulami2023}}
\label{tab:parameters}
\begin{tabular}{c c p{0.6\textwidth}}
\toprule
Parameter & Value & Description \\
\midrule
\( a \) & 1 & Coefficient for radiative heat loss efficiency. \\
\( b \) & 0.3 & Magnitude of incoming solar radiation flux. \\
\( s_0 \) & 0.1 & Surface albedo at high phytoplankton density. \\
\( s_1 \) & 0.95 & Surface albedo in the absence of phytoplankton. \\
\( \alpha_1 \) & 3 & Sensitivity of albedo to biomass fluctuations. \\
\( T_0 \) & 1 & Low-temperature growth limitation. \\
\( \alpha_2 \) & 1 & High-temperature growth suppression. \\
\( \gamma \) & 1 & Thermal inertia of the upper water column. \\
\( \mu \) & 0.1 & Phytoplankton mortality and loss rate. \\
\( b_1 \) & 2.1 & Reference intrinsic growth rate. \\
\bottomrule
\end{tabular}
\end{table}

%\rev{ 
The parameter values listed in Table~\ref{tab:parameters} are adopted from the temperature–phytoplankton
model in the context of Arctic UIB
dynamics. These values correspond to a biologically realistic regime in which
light-limited phytoplankton growth, nutrient feedback, and temperature-dependent
mortality interact to produce bistability between a low-biomass background
state and a bloom state.

The suitability of the model for UIB dynamics is mainly encoded in the biomass-dependent radiative feedback term $S(u)$ and in the chosen bistable parameter regime. In particular, $s_1=0.95$ represents an effective high-albedo, strongly light-limited state in the absence of phytoplankton, while $s_0=0.1$ corresponds to enhanced absorption at high biomass. Thus, as $u$ increases, the decrease of $S(u)$ strengthens shortwave absorption through the factor $1-S(u)$, improving the effective thermal/light conditions for growth. Together with the temperature-dependent growth function $g(T)$ and the loss parameter $\mu$, this feedback generates coexistence of a low-biomass under-ice state and a bloom state, making the model appropriate for investigating abrupt UIB onset.

In ecological terms, the low-biomass equilibrium corresponds to a background
under-ice state in which light limitation and mortality dominate growth,
whereas the high-biomass equilibrium represents a bloom regime sustained by
enhanced absorption of shortwave radiation and positive growth feedback.
The function $S(u)$ introduces a biomass–radiation feedback: as $u$ increases,
light absorption in the upper ocean is enhanced, which in turn modifies the
effective growth conditions. This feedback mechanism enables the coexistence
of two stable equilibria separated by a saddle, providing a deterministic
mechanism for abrupt bloom onset. 

In the absence of stochastic perturbations, trajectories remain confined to the basin of their initial state, while environmental variability introduced in Section~2.2 allows transitions between these regimes. 
In the present work, we do not attempt to re-estimate the ecological parameters. Instead, we adopt the set of parameters used in the previously validated temperature–phytoplankton model~\cite{Alsulami2023}. This allows us to focus on the geometric
structure of stochastic transitions rather than on parameter identification. The chosen control-parameter range for $b_1$ lies entirely within the bistable window identified in the deterministic system, ensuring the existence of two stable equilibria separated by a saddle equilibrium and, therefore, allowing for a well-defined stochastic separatrix analysis.

\subsection{Stochastic model}\label{subsec:stoch_model}
%high-frequency stochastic fluctuations
To account for high-frequency stochastic fluctuations such as rapid radiation changes, mixed-layer disturbances, and nutrient pulses, we extend the deterministic system in Eqs.~\eqref{eq:det_T}--\eqref{eq:det_u} to the following two-dimensional SDE~\cite{Wouters2016,Franzke2015StochClimate}:
\begin{align}
\mathrm{d}T_t &= F_T(T_t,u_t)\,\mathrm{d}t + \sigma\,\mathrm{d}W^{(T)}_t, \label{eq:sde_T}\\
\mathrm{d}u_t &= F_u(T_t,u_t)\,\mathrm{d}t + \sigma\sqrt{u_t+\delta}\,\mathrm{d}W^{(u)}_t. \label{eq:sde_u}
\end{align}
where $W_t^{(T)}$ and $W_t^{(u)}$ denote independent standard Wiener processes, and $\sigma$ represents the noise intensity. The stochastic differential
equations are interpreted in the It\^o sense, which leads to the
Fokker--Planck equation (FPE) for the probability density derived in Section~\ref{subsec:fpe}~\cite{Risken1996}.

%\revgreen{In the deterministic limit ($\sigma = 0$), the equilibrium structure
%shown by the black curves in Fig.~1 exhibits a classical saddle–node
%bifurcation scenario with respect to the control parameter $b_1$.
%Two stable equilibria (a low-biomass background state $E_1$ and
%a high-biomass bloom state $E_3$) coexist within the interval
%$b_1 \in [1.996,\, 2.471]$, separated by an unstable saddle $E_2$.
%At the endpoints of this interval, the positive equilibrium is
%created or annihilated through saddle–node bifurcations.
%Thus, bistability occurs only for sufficiently overcritical values
%of $b_1$, and outside this interval the system admits a single
%stable steady state.}

The temperature component uses additive noise to represent fluctuations in external thermal forcing. The biomass component uses the noise term $\sqrt{u_t+\delta}$, which matches the square root scaling in demographic stochasticity and reduces the probability that discrete time simulations cross $u=0$~\cite{Risken1996}.

The term $\sqrt{u+\delta}$ is selected over $\sqrt{u}+\delta$ to ensure that regularisation remains localised to a small neighbourhood of $u\approx 0$. While the latter would introduce a persistent additive diffusion across the entire domain, the former ensures the regularisation effect is negligible outside the immediate vicinity of the origin. Setting the regularisation parameter $\delta = 10^{-4}$ ensures that $\delta \ll u_{\min}$, maintaining a negligible impact on the dynamics for $u > 0$. Both components share the same noise intensity $\sigma$ to facilitate parameter analysis; supplementary computations with $\sigma_T \neq \sigma_u$ indicate that the qualitative conclusions remain unchanged (Section~\ref{subsec:robustness}).

\subsection{Fokker--Planck equation}\label{subsec:fpe}
Let $p=p(T,u,t)$ be the joint probability density of the process $(T_t,u_t)$. The FPE associated with the SDE system Eqs.~\eqref{eq:sde_T}--\eqref{eq:sde_u} is given by:
\begin{equation}
\frac{\partial p}{\partial t}
= \mathcal{L}^* p
= -\frac{\partial}{\partial T}\big(F_T\,p\big)
  -\frac{\partial}{\partial u}\big(F_u\,p\big)
  +\frac{\sigma^2}{2}\frac{\partial^2 p}{\partial T^2}
  +\frac{1}{2}\frac{\partial^2}{\partial u^2}\Big(\sigma^2(u+\delta)\,p\Big).
\label{eq:fpe_2d}
\end{equation}
where $\mathcal{L}^*$ is the forward Kolmogorov generator.
The second–order (diffusion) terms in Eq.~\eqref{eq:fpe_2d}
correspond to the diffusion tensor associated with the It\^o noise coefficient matrix
\[
G(T,u)=\mathrm{diag}(\sigma,\sigma\sqrt{u+\delta}).
\]
The resulting diffusion tensor is
\begin{equation}
a(T,u)=G G^\top
=\mathrm{diag}(\sigma^2,\sigma^2(u+\delta)).
\label{eq:diff_tensor}
\end{equation}

We impose reflecting zero-flux boundary conditions on a bounded domain $\Omega$:
\begin{equation}
n\cdot J=0 \qquad \text{on }\partial\Omega,
\label{eq:no_flux_bc}
\end{equation}
where $J$ denotes the probability flux and $n$ is the outward unit normal. For subsequent backward problems, the boundary $\partial\Omega$ is assigned homogeneous Neumann conditions. 

The operator $\mathcal{L}$ is defined by its action on a smooth test
function $\varphi(T,u)$:
\begin{equation}
\mathcal{L}\varphi
= F_T\,\partial_T \varphi + F_u\,\partial_u \varphi
+ \frac{\sigma^2}{2}\,\partial_{TT}\varphi
+ \frac{\sigma^2}{2}(u+\delta)\,\partial_{uu}\varphi .
\label{eq:generator_L}
\end{equation}
This operator $\mathcal{L}$ is the backward Kolmogorov generator associated with the stochastic differential system in Eqs.(4)--(5). It acts on observables, whereas its adjoint $\mathcal{L}^*$ acts on probability densities, satisfying the adjoint relation~\cite{Risken1996}
\[
\int_{\Omega} \varphi(T,u)\,(\mathcal{L}^\ast p)(T,u,t)\,dT\,du
=
\int_{\Omega} (\mathcal{L}\varphi)(T,u)\,p(T,u,t)\,dT\,du .
\]

The committor function and MFPT will be formulated as boundary value problems governed by this generator $\mathcal{L}$ in Sections~\ref{subsec:committor_separatrix} and Section~\ref{subsec:mfpt_def}. 
\subsection{Stochastic bifurcation}\label{subsec:sbd_expect}
\begin{figure}[t]
  \centering
  \includegraphics[width=0.8\linewidth]{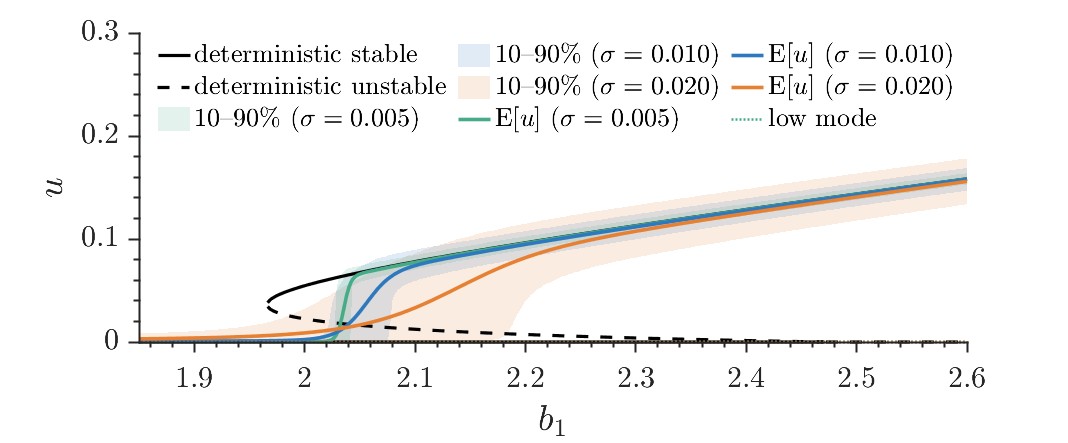}
  \caption{Expectation based stochastic bifurcation diagram.
Black curves represent the deterministic equilibrium branches: solid lines denote stable equilibria ($E_1$ and $E_3$), while the dashed line denotes the unstable saddle ($E_2$). Coloured solid curves depict the stationary expectation $\bar{u}(b_1,\sigma)$ for noise intensities $\sigma\in\{0.005, 0.010, 0.020\}$. Shaded bands indicate the corresponding $10\%$--$90\%$ quantile ranges of the stationary distribution.}
  \label{fig:sbd_expect}
\end{figure}

To examine how stationary statistics evolve with the control parameter
under stochastic perturbations, we construct a stochastic bifurcation
diagram, shown in Fig.~\ref{fig:sbd_expect}, based on the expectation of the
stationary distribution $\bar{u}(b_1,\sigma)$.

Under a quasi-steady reduction for temperature, we set $\dot T=0$ in the deterministic $T$-equation to obtain
\begin{equation}
\label{eq:Tstar}
T^\ast(u)=\left[\frac{b}{a}\bigl(1-S(u)\bigr)\right]^{1/4}.
\end{equation}
Substituting $T^\ast(u)$ into $g(T)$ and then into the $u$-equation yields
a closed one-dimensional It\^o SDE for $u$, which retains the multiplicative
noise structure $\sigma\sqrt{u+\delta}\,dW_t$ from the original two-dimensional system.

%Substituting $T^\ast(u)$ into $g(T)$ and then into the $u$-equation yields a closed 1D It\^o SDE
\begin{equation}
\label{eq:1d_sde}
du = f(u;b_1)\,dt + \sigma\sqrt{u+\delta}\,dW_t,
\qquad
f(u;b_1)\equiv u\,[g(T^\ast(u))-u]-\mu u.
\end{equation}
%Under this reduction, the $T$-component is eliminated and the
%resulting one-dimensional It\^o SDE inherits the multiplicative
%noise structure from the $u$-equation.

Let the effective one-dimensional FPE diffusion coefficient be
\begin{equation}
\label{eq:D_def}
D(u)=\tfrac{\sigma^2}{2}(u+\delta).
\end{equation}
The coefficient $D(u)$ is obtained from the $uu$-component
of the diffusion tensor $a(T,u)$ defined in Eq.~(7),
namely $a_{uu}(T,u)=\sigma^2(u+\delta)$,
so that $D(u)=\tfrac12 a_{uu}(T^*(u),u)$.

In the one-dimensional reduction, the FPE operator
takes the standard form
\[
\partial_t p = -\partial_u(fp)
+ \partial_u^2(D(u)p)\]

On a bounded interval with reflecting (zero-flux) boundaries,
the stationary density of the one-dimensional reduction
is obtained by setting $\partial_t p = 0$ in the corresponding
one-dimensional FPE derived from Eq.~\eqref{eq:fpe_2d},
which yields
\begin{equation*}
0=-\partial_u J(u), \qquad J(u)=f(u;b_1)\,p(u)-\partial_u\!\big(D(u)\,p(u)\big).
\end{equation*}
The reflecting boundaries enforce $J\equiv 0$, hence
\begin{equation*}
\partial_u\!\big(D p\big)=f\,p
\;\;\Longrightarrow\;\;
\frac{d}{du}\ln\!\big(D p\big)=\frac{f}{D}.
\end{equation*}
Integrating gives
\begin{equation}
\label{eq:1d_stationary_p}
p(u)=\frac{C}{D(u)}\,
\exp\!\left(\int^{u}\frac{f(z;b_1)}{D(z)}\,dz\right)
= \frac{C'}{u+\delta}\,
\exp\!\left(\int^{u}\frac{2 f(z;b_1)}{\sigma^2(z+\delta)}\,dz\right)\end{equation}
where $C$ (equivalently $C'$) normalises $\int p(u)\,du=1$.
We construct an expectation-centred stochastic bifurcation diagram via
\begin{equation}
\label{eq:ubar_def}
\bar u(b_1,\sigma)=\int u\,p(u)\,du
\end{equation}
and depict dispersion with the $10\%$ and $90\%$ quantiles. 

This expectation-based representation is robust to mode switching and avoids spurious discontinuities from tracking the maximum of a bimodal probability density function~\cite{Arnold1998Bifurcation}.
%%%%%%%%%%%[new]%%%%%
In the deterministic limit ($\sigma=0$), the system exhibits bistability within the parameter window $b_1 \in [1.996, 2.471]$, bounded by saddle-node bifurcations at both critical ends where the unstable saddle $E_2$ coalesces with either the stable background state $E_1$ or the bloom state $E_3$. Within this bistable interval, for example at the representative value $b_1 = 2.1$, the coexisting stable equilibria are located at $E_1 = (0.350, 0.000)$ and $E_3 \approx (0.511, 0.078)$, separated by the saddle at $E_2 \approx (0.395, 0.012)$.
%%In the deterministic limit ($\sigma=0$), the system exhibits a bistable structure for $b_1 \in [1.996,\,2.471]$, where the stable background state $E_1$ and bloom state $E_3$ are separated by an unstable saddle $E_2$. At the representative value $b_1 = 2.1$, the coexisting stable equilibria are located at $E_1 = (0.350, 0.000)$ and $E_3 \approx (0.511, 0.078)$, with the saddle at $E_2 \approx (0.395, 0.012)$.

%Under stochastic perturbations, the transition of $\bar{u}$ with $b_1$ becomes a continuous shift rather than a discontinuous jump. As noise intensity $\sigma$ increases, the quantile bands widen, reflecting increased stationary distribution dispersion, and the transition interval shifts toward larger $b_1$.

%\revgreen{Under stochastic perturbations, the transition of $\bar u$ with $b_1$
%becomes a continuous shift rather than a discontinuous jump.
%As noise intensity $\sigma$ increases, the quantile bands widen,
%reflecting increased stationary dispersion, and the transition
%interval shifts toward larger $b_1$.
%Noise therefore regularises the observable bifurcation diagram
%by smoothing expectation-based statistics, while the deterministic
%saddle--node structure persists in the limit $\sigma \to 0$.
%For small but finite noise, bistability is replaced by metastability:
%the stationary distribution becomes bimodal and transitions between
%states occur with finite probability.
%The smooth variation of $\bar u(b_1,\sigma)$ thus reflects changes
%in the relative occupation probabilities of the two metastable states,
%rather than the disappearance of tipping.}

Under stochastic perturbations, the transition of $\bar{u}$ with $b_1$ becomes a continuous shift rather than a discontinuous jump. As noise intensity $\sigma$ increases, the quantile bands widen, reflecting increased stationary dispersion, and the transition interval shifts toward larger $b_1$. Noise therefore regularises the observable bifurcation diagram by smoothing expectation-based statistics, while the deterministic saddle–node structure persists in the limit $\sigma \to 0$. From the probabilistic viewpoint, the stochastic bifurcation manifests itself through a qualitative change in the stationary probability distribution. In particular, as the control parameter $b_1$ varies, the stationary density transitions from a unimodal distribution associated with a single dominant state to a bimodal distribution corresponding to two metastable states. For small but finite noise, bistability is replaced by metastability: the stationary distribution becomes bimodal and transitions between states occur with finite probability. The smooth variation of $\bar{u}(b_1,\sigma)$ therefore reflects continuous changes in the relative occupation probabilities of the two metastable states, rather than the disappearance of tipping.

\subsection{Stationary marginal distributions}\label{subsec:stat_density}

\begin{figure}[H]
    \centering
    \includegraphics[width=0.95\textwidth]{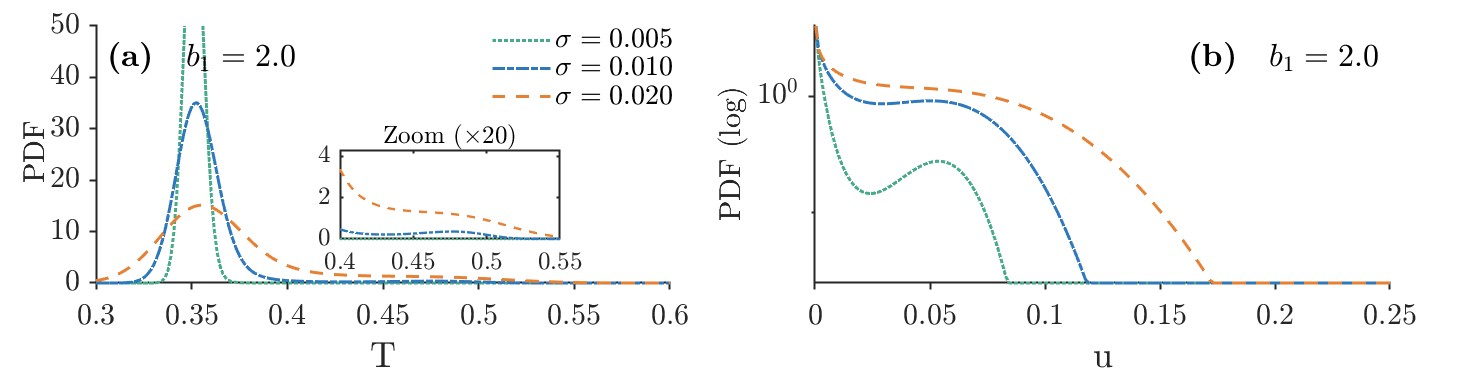} \\
   % \vspace{0.2cm}
    \includegraphics[width=0.95\textwidth]{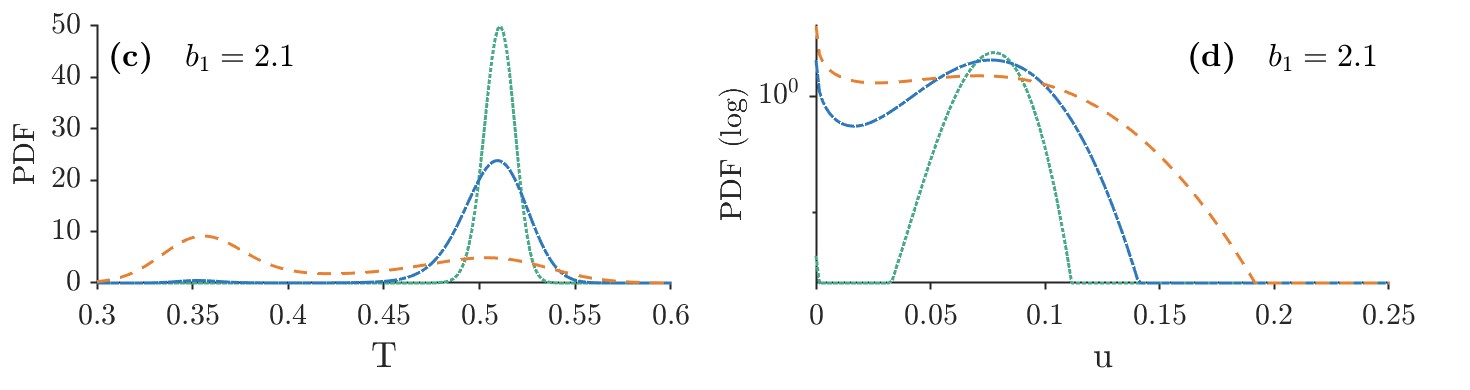} \\
  %  \vspace{0.2cm}
    \includegraphics[width=0.95\textwidth]{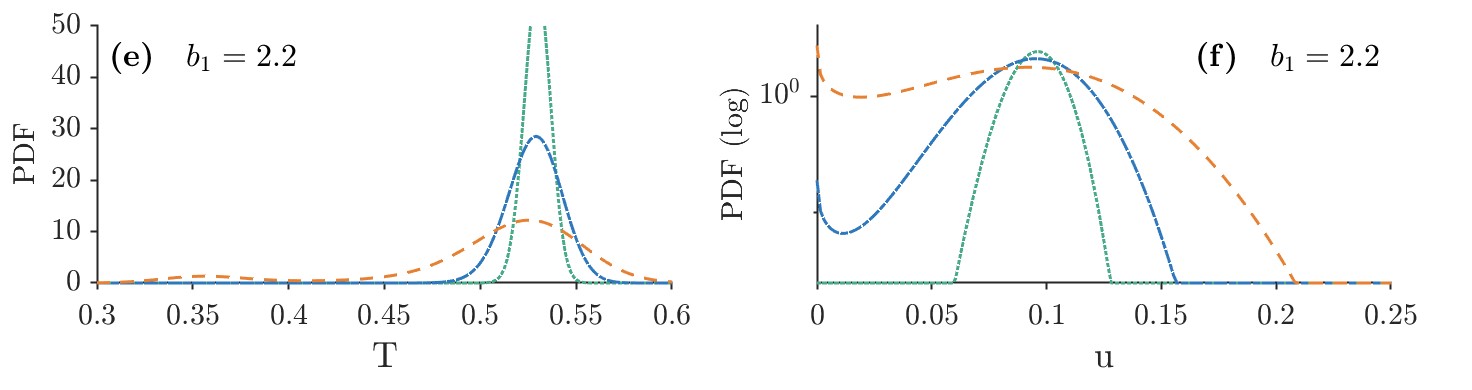}
  \caption{Approximate stationary marginal probability density functions (PDFs). 
Left column (a, c, e): Temperature marginals $p_T(T)$, where the inset in (a) shows a 20-fold vertical magnification of the region $T \in [0.4, 0.55]$ to resolve the weak $E_3$ peak. 
Right column (b, d, f): Biomass marginals $p_u(u)$ on a logarithmic vertical axis. 
Rows from top to bottom correspond to control parameter values $b_1=2.0$, $b_1=2.1$, and $b_1=2.2$. 
Curves represent noise intensities $\sigma=0.005$ (green dotted), $0.010$ (blue dash-dotted), and $0.020$ (orange dashed).}
    \label{fig:stat_marginals}
\end{figure}

To elucidate the microscopic origin of the stochastic bifurcation diagram trends, we examine the stationary
marginals of the joint stationary density $p(T,u)$ on $\Omega=[T_{\min},T_{\max}]\times[u_{\min},u_{\max}]$:
\begin{equation}
\label{eq:pt_pu}
p_T(T)=\int_{u_{\min}}^{u_{\max}} p(T,u)\,du,\qquad
p_u(u)=\int_{T_{\min}}^{T_{\max}} p(T,u)\,dT.
\end{equation}
and display them in Fig.~\ref{fig:stat_marginals}.
These marginals are computed either from long SDE simulations with reflecting
boundaries or from the stationary FPE with zero-flux
boundary conditions; the two approaches agree within numerical tolerance. The biomass marginal $p_u(u)$ is plotted on a logarithmic vertical axis to resolve the probability mass accumulating near the extinction boundary $u \approx 0$.

As noise intensity $\sigma$ increases, the probability peaks broaden and their heights decrease, enhancing density near $E_1$; while larger $b_1$ values shift the distribution modes toward $E_3$. Throughout the bistable regime, the double peaks structure is most pronounced at $b_1=2.1$ and less distinct at $b_1=2.0$ and $b_1=2.2$, with the valley corresponding to the transition region. For $b_1 = 2.0$, although the deterministic system lies within the bistable interval, the secondary temperature peak is considerably weaker due to asymmetry between the potential wells. This effect is consistent with the general observation that in the zero noise limit, unless symmetries are present in the system, all the probability will concentrate around only one peak and the population of the secondary minima of the quasi-potential becomes negligible~\cite{LucariniBodai2020}. Furthermore, because the temperature coordinates of the two equilibria are closer than their biomass coordinates, marginalisation over $u$ additionally reduces the visibility of the second mode in $p_T(T)$. The inset in Fig.~\ref{fig:stat_marginals}(a) reveals this weak secondary peak, confirming that bimodality persists at the joint level.

\section{Geometric indicator}\label{sec:methodology}

\subsection{Committor function and stochastic separatrix}\label{subsec:committor_separatrix}
In the bistable regime, stochastic perturbations can induce noise-induced transitions between the background state and the bloom state. To quantify the likelihood of a trajectory reaching a specific attractor first, we introduce the forward committor function~\cite{WeinanE2006}.

Let $X_t=(T_t,u_t)$ solve the system Eqs.~\eqref{eq:sde_T}--\eqref{eq:sde_u}. Within the bounded computational domain $\Omega \subset \mathbb{R}^2$, we define two disjoint closed neighbourhoods $R_{E_1}$ and $R_{E_3}$ containing the stable equilibria $E_1$ and $E_3$, respectively. These neighbourhoods are specified as axis-aligned ellipses centred at the corresponding equilibria in our SDE model. For any initial condition $x \in \Omega$, the first hitting times are defined as

\[
\tau_{R_{E_1}} = \inf\{t \geq 0 : X_t \in R_{E_1}\}, \qquad
\tau_{R_{E_3}} = \inf\{t \geq 0 : X_t \in R_{E_3}\}.
\]
The forward committor function $q(T,u; b_1,\sigma)$ is the probability that a trajectory starting from $(T,u)$ hits $R_{E_3}$ before $R_{E_1}$:
\begin{equation}
q(T,u; b_1,\sigma) = \mathbb{P}\!\left(\tau_{R_{E_3}} < \tau_{R_{E_1}} \,\big|\, X_0 = (T,u)\right).
\label{eq:committor_def}
\end{equation}

By Dynkin's formula, $q$ satisfies the backward Kolmogorov elliptic boundary value problem~\cite{Risken1996}:
\begin{equation}
\left\{
\begin{aligned}
\mathcal{L} q &= 0, && (T,u) \in \Omega \setminus (R_{E_1} \cup R_{E_3}), \\
q|_{R_{E_1}} &= 0, \quad q|_{R_{E_3}} = 1, \\
\partial_n q|_{\partial\Omega} &= 0,
\end{aligned}
\right.
\label{eq:committor_bvp}
\end{equation}
where $\mathcal{L}$ is the generator given by Eq.~\eqref{eq:generator_L}, and homogeneous Neumann conditions are imposed on the outer boundary $\partial\Omega$. 

Values of $q$ close to $0$ indicate a higher probability of first reaching $R_{E_1}$, whereas values close to $1$ signify a higher probability of reaching $R_{E_3}$. We define the stochastic separatrix $\Gamma$ as the isocommittor $q=1/2$:

\begin{equation}
\Gamma(b_1,\sigma) = \{ (T,u) \in \Omega : q(T,u; b_1,\sigma) = 1/2 \}.
\label{eq:Gamma_def}
\end{equation}

Along $\Gamma$, the eventual outcome is maximally uncertain~\cite{AntoniouSchwartz2009JCP}. Unlike the deterministic separatrix, $\Gamma$ is defined at finite noise intensity and does not necessarily pass through the saddle $E_2$. Numerical solution of Eq.~\eqref{eq:committor_bvp} follows the setup in~\ref{app:setup}.

%\subsection{Geometric indicator}\label{subsec:ews_geom}
\subsection{Geometric indicator $EWS_{\mathrm{geom}}$}\label{subsec:ews_geom}
The stochastic separatrix $\Gamma$ forms the centre of a probabilistic transition layer rather than a sharp boundary. We define a geometric indicator $EWS_{\mathrm{geom}}$ that measures the spatial width of this layer using the normal gradient of $q$, where an increase in this width reflects the weakening of basin stability and a higher likelihood of noise-induced transitions.
%%%%%%%%%%%above new [where an increase in this width reflects the weakening of basin stability and a higher likelihood of noise-induced transitions.]

For any $x \in \Gamma$, the unit normal vector pointing toward the bloom state is defined as 
\[
n(x) = \frac{\nabla q(x)}{\|\nabla q(x)\|}.
\]

Let $\xi$ denote displacement along $n(x)$. A first-order Taylor expansion yields:
\begin{equation}
q(x + \xi n) \approx q(x) + \xi \, \nabla q(x) \cdot n(x) = \frac{1}{2} + \xi \, \|\nabla q(x)\|.
\label{eq:q_normal_expansion}
\end{equation}
\begin{figure}[t]
\centering
\includegraphics[width=0.7\textwidth]{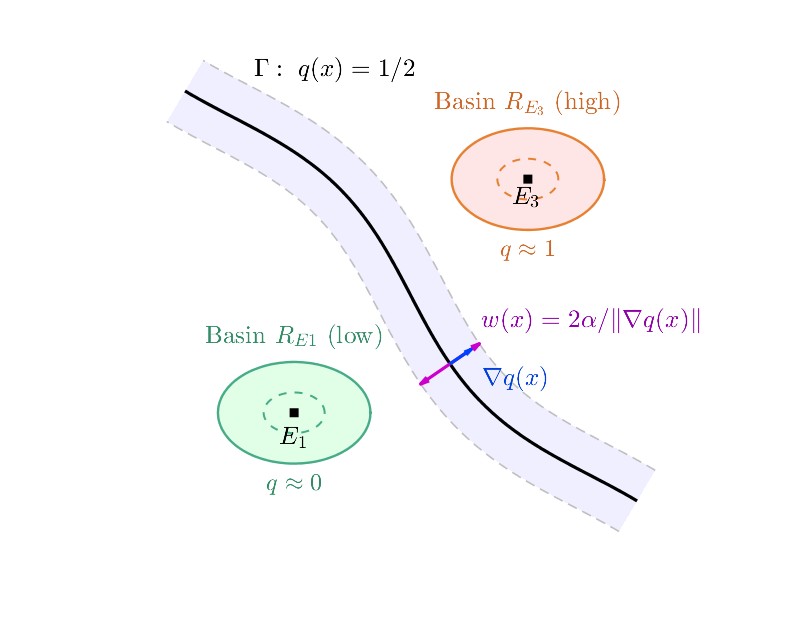}
\caption{Schematic representation of the stochastic separatrix and its associated transition layer. The stochastic separatrix $\Gamma = \{q = 1/2\}$ (black curve) separates the background basin ($q \approx 0$) from the bloom basin ($q \approx 1$). The transition layer (shaded region) corresponds to $q \in [0.4,0.6]$ with $\alpha = 0.1$. The local width $w_\alpha(x)$ is defined as the normal distance between the $q=0.4$ and $q=0.6$ boundaries. The gradient $\nabla q(x)$ points toward the bloom state. The geometric early warning indicator $\mathrm{EWS}_{\mathrm{geom}}$ is the arc-length average of this local width along $\Gamma$.}
%Conceptual illustration of the geometric indicator...
\label{fig:sep_concept}
\end{figure}
Introducing a probability half-width $\alpha \in (0,1/2)$, the distances to the boundaries $q = 1/2 \pm \alpha$ satisfy
\begin{equation}
\xi_+(x;\alpha) \approx \xi_-(x;\alpha) \approx \frac{\alpha}{\|\nabla q(x)\|}.
\label{eq:xi_pm}
\end{equation}

The local full width $w_\alpha(x)$ is thus defined as:
\begin{equation}
w_\alpha(x) := \xi_+(x;\alpha) + \xi_-(x;\alpha) \approx \frac{2\alpha}{\|\nabla q(x)\|}.
\label{eq:local_width}
\end{equation}

Thus, smaller $\|\nabla q\|$ corresponds to slower spatial variation in probability and a wider transition layer.
The quantity $w_{\alpha}(x)$ represents a local geometric measure of the transition-layer thickness at a point on the stochastic separatrix. The early warning indicator is not this local width itself, but the global scalar obtained by averaging $w_{\alpha}$ along $\Gamma$, as defined below.

To obtain a global scalar measure of basin-boundary sharpness across the control parameter space $(b_1,\sigma)$, we define the geometric early warning indicator $\mathrm{EWS}_{\mathrm{geom}}$ as the arc-length average of the local transition-layer width $w_{\alpha}$ along $\Gamma$:

%To compare layer thickness across the control parameter space $(b_1,\sigma)$, we define the  geometric indicator $EWS_{\mathrm{geom}}$ as the arc-length average of $w_\alpha$ along $\Gamma$:
\begin{equation}
EWS_{\mathrm{geom}}(b_1,\sigma;\alpha)
:= \frac{1}{L(\Gamma)} \int_{\Gamma} w_\alpha(s) \, \mathrm{d}s
\approx \frac{1}{L(\Gamma)} \int_{\Gamma} \frac{2\alpha}{\|\nabla q(s)\|} \, \mathrm{d}s.
\label{eq:ews_geom_def}
\end{equation}
Since $EWS_{\mathrm{geom}}$ scales linearly with $\alpha$, we fix $\alpha = 0.1$ corresponding to the core band $q \in [0.4,0.6]$. Gradient regularisation and implementation details are provided in~\ref{app:ews}. An increase in $\mathrm{EWS}_{\mathrm{geom}}$ corresponds to a broader probabilistic transition layer and therefore indicates reduced basin stability and increased susceptibility to noise-induced switching.

Furthermore, although curvature of the stochastic separatrix induces asymmetry between the two half-widths, this effect does not alter the dominant $O(\alpha)$ scaling of the total width $w_\alpha$; the complete second-order analysis is provided in~\ref{app:second_order}. 

\subsection{Separatrix shift}\label{subsec:geom_baselines}
To quantify the overall shift of the stochastic separatrix $\Gamma$, we introduce two baseline quantities based on mean distances, namely the mean distance to the deterministic boundary ($MDB$) and the mean distance to the stochastic separatrix ($MDS$). These metrics measure the deviation of $\Gamma$ relative to the deterministic limit $\Gamma_{\mathrm{det}}(b_1)$ and to a reference background stochastic state $\Gamma(b_1,\sigma_{\mathrm{ref}})$, respectively.

For any two curves $\Gamma_1$ and $\Gamma_2$ in phase space, we define the
arc-length averaged minimal distance from $\Gamma_1$ to $\Gamma_2$ as
\begin{equation}
D(\Gamma_1 \to \Gamma_2)
= \frac{1}{L(\Gamma_1)} \int_{\Gamma_1} \operatorname{dist}\!\left(x,\Gamma_2\right) \, \mathrm{d}s,
\label{eq:arc_distance}
\end{equation}
where $x\in\Gamma_1$ denotes the integration point and
\[
\operatorname{dist}(x,\Gamma_2) := \inf_{y\in\Gamma_2}\|x-y\|_2
\]
is the Euclidean distance from $x$ to the closest point on $\Gamma_2$. Note that $D(\Gamma_1\to\Gamma_2)$ is a directed distance and generally $D(\Gamma_1\to\Gamma_2)\neq D(\Gamma_2\to\Gamma_1)$.

The shift indicators are then given by:
\begin{equation}
MDB(b_1,\sigma) = D\!\left(\Gamma(b_1,\sigma) \to \Gamma_{\mathrm{det}}(b_1)\right),
\label{eq:MDB_def}
\end{equation}
\begin{equation}
MDS(b_1,\sigma) = D\!\left(\Gamma(b_1,\sigma) \to \Gamma(b_1,\sigma_{\mathrm{ref}})\right), \quad \sigma_{\mathrm{ref}} = 0.005.
\label{eq:MDS_def}
\end{equation}

%%%%%%%%%%%%%%%%%%%%%%%%%%
\section{Geometric analysis}\label{sec:results}
\subsection{Separatrix geometry}\label{subsec:sep_geometry}
\begin{figure}[h]
\centering
\includegraphics[width=\linewidth]{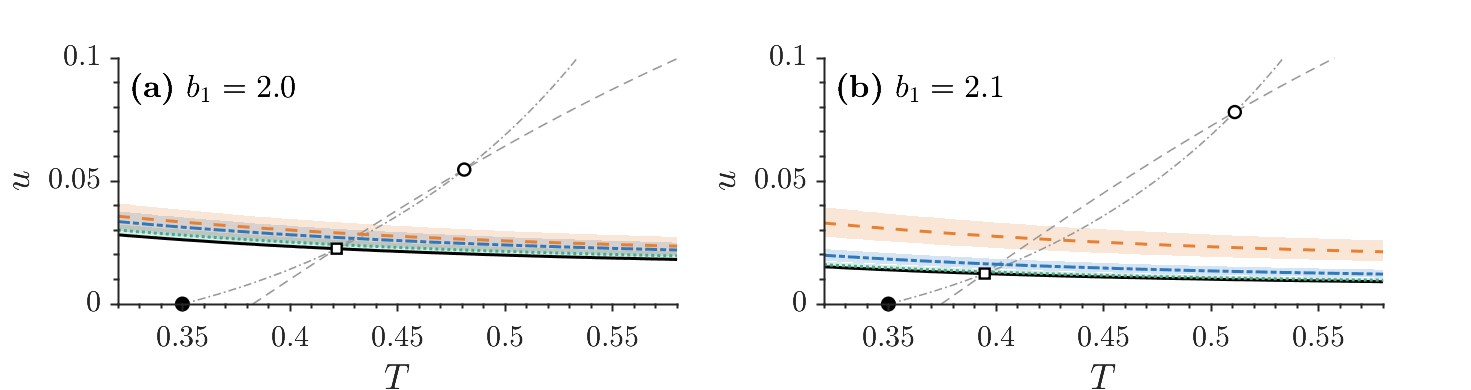}
    \centering
    \includegraphics[width=\linewidth]{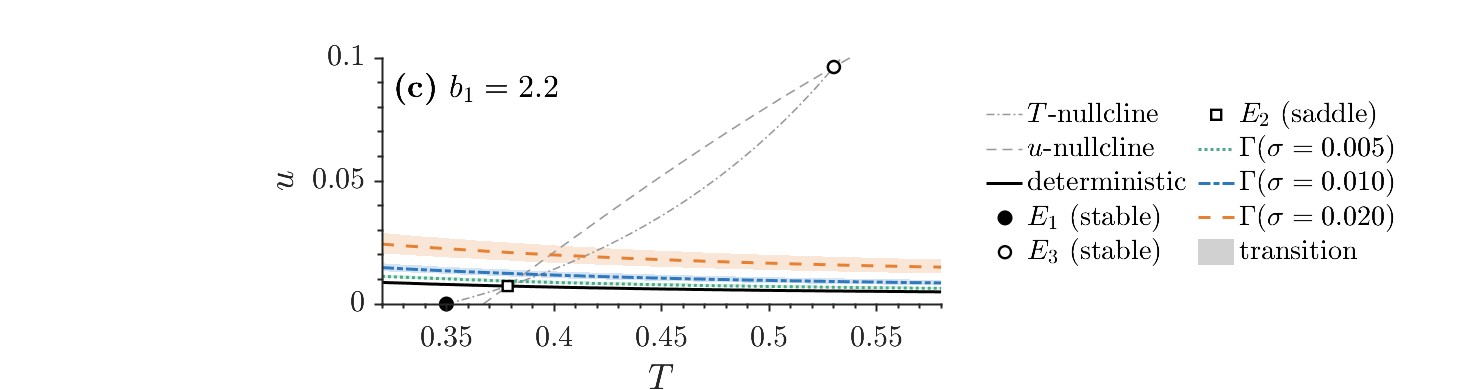}
\caption{
Deterministic and stochastic separatrices in the $(T,u)$ phase plane for (a) $b_1=2.0$, (b) $b_1=2.1$, and (c) $b_1=2.2$. 
The deterministic separatrix is shown as a solid black curve, $T$- and $u$-nullclines are indicated by gray lines.
Coloured curves denote the stochastic separatrices $\Gamma(b_1,\sigma)$ for $\sigma=0.005$ (green dotted), $\sigma=0.010$ (blue dashed--dotted), and $\sigma=0.020$ (orange dash). 
The transition layer ($q\in[0.4,0.6]$) is shown as shaded regions using the same colour scheme. 
Symbols $E_1$ (solid circle), $E_3$ (open circle), and $E_2$ (open square) represent stable equilibria and the saddle point, respectively. }
\label{fig:sep_overlay}
\end{figure}

We analyse the geometric relationship between the position of the deterministic separatrix $\Gamma_{\mathrm{det}}$ and the stochastic separatrix $\Gamma$ in the $(T,u)$ phase plane. The Fig.~\ref{fig:sep_overlay} illustrates $\Gamma_{\mathrm{det}}$ together with $\Gamma$ under three noise intensities, with shaded regions indicating the transition layer.

Within the investigated parameter range, $\Gamma$ and $\Gamma_{\mathrm{det}}$
maintain similar geometric shapes. 
In the weak-noise limit $\sigma \to 0$, the stochastic separatrix converges to the deterministic separatrix, consistent with the deterministic limit of the stochastic dynamics.
%In the weak-noise limit $\sigma \to 0$, the stochastic separatrix converges to the deterministic separatrix, so that the distance between them vanishes. 
For fixed $b_1$,
increasing $\sigma$ shifts $\Gamma$ toward larger $u$, reflecting
a noise-induced bias toward $E_1$. This shift implies elevated extinction
risk for phytoplankton under stronger environmental fluctuations, as
the basin of attraction for the bloom state contracts. Simultaneously,
the transition layer broadens, indicating an expanded region where basin
transitions occur with comparable likelihood.

%Within the investigated parameter range, $\Gamma$ and $\Gamma_{\mathrm{det}}$ maintain similar geometric shapes. For fixed $b_1$, increasing $\sigma$ shifts $\Gamma$ toward larger $u$, reflecting a noise-induced bias toward $E_1$. This shift implies elevated extinction risk for phytoplankton under stronger environmental fluctuations, as the basin of attraction for the bloom state contracts. Simultaneously, the transition layer broadens, indicating an expanded region where basin transitions occur with comparable likelihood.

\subsection{Separatrix measures}\label{subsec:sep_shift}

\begin{figure}[t]
\centering
\includegraphics[width=\linewidth]{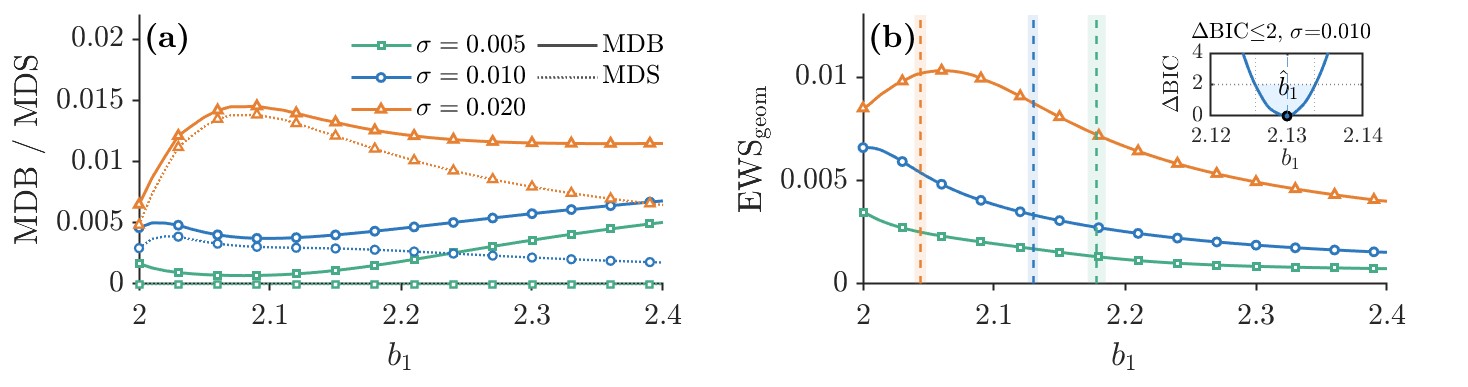}
\caption{Geometric measures under varying noise intensity. (a) Positional shift measures $MDB$ (solid lines) and $MDS$ (dotted lines) as functions of $b_1$. (b) $EWS_{\mathrm{geom}}$ versus $b_1$, showing the BIC-selected optimal breakpoints $\hat b_1$ (vertical dashed lines) and the corresponding warning intervals (shaded bands, $\Delta \mathrm{BIC} \le 2$). Curves correspond to $\sigma=0.005$ (green), $\sigma=0.010$ (blue), and $\sigma=0.020$ (orange). Inset: BIC profile for $\sigma = 0.010$.
}
\label{fig:sep_shift}
\end{figure}

We quantify the geometric variations shown in Fig.~\ref{fig:sep_shift}
using the positional shift measures $MDB$ and $MDS$
(Section~\ref{subsec:geom_baselines})
and the transition-layer width measure $\mathrm{EWS}_{\mathrm{geom}}$
(Section~\ref{subsec:ews_geom}).

The Fig.~\ref{fig:sep_shift}(a) shows that the noise-induced shift is non-uniform along the $b_1$ axis. At $b_1 = 2.08$, the MDB expands by a factor of $21.0$ as $\sigma$ increases from $0.005$ to $0.020$, while at $b_1 = 2.40$ it increases by a factor of only $2.28$. The ratio $MDS/MDB$ decays from $0.95$ to $0.56$ across the $b_1$ range, indicating a crossover from a noise-dominated regime to background regime $\sigma_{\mathrm{ref}}$.

The observed displacement of the stochastic separatrix relative to the
deterministic one is partly associated with the multiplicative structure
of the biomass noise term $\sigma\sqrt{u+\delta}$, which makes the
effective diffusion strength state dependent. This asymmetry modifies
the balance of stochastic fluxes near the deterministic separatrix and
leads to a systematic shift of the stochastic basin boundary.
However, the phenomenon is not solely a consequence of multiplicative
noise. Additional computations with alternative noise structures
(e.g., independent noise intensities in the temperature and biomass
components) indicate that the qualitative behaviour of the separatrix
shift and the associated geometric indicators remains unchanged.
Thus, the distance measures $MDB$ and $MDS$ capture a general
stochastic deformation of the basin boundary rather than an artifact
specific to the chosen noise parametrisation.

The Fig.~\ref{fig:sep_shift}(b) presents that $\mathrm{EWS}_{\mathrm{geom}}$ decreases with $b_1$ and increases with $\sigma$, consistent with the transition layer broadening observed in Fig.~\ref{fig:sep_overlay}.  At $b_1 = 2.40$, $MDB$ is small while $\mathrm{EWS}_{\mathrm{geom}}$ is broaden, indicating that the positional shift and the transition layer width are not strictly coupled. 

Furthermore, we apply a two-segment linear fit to the $\mathrm{EWS}_{\mathrm{geom}}$ curves. Using the Bayesian Information Criterion (BIC), we identify the optimal breakpoint $\hat b_1$, which corresponds to the abrupt change in the rate of increase of $\mathrm{EWS}_{\mathrm{geom}}$, and determine the early warning interval based on the criterion $\Delta\mathrm{BIC}\le 2$ (specific values are listed in~\ref{app:bic}). Results indicate that as $\sigma$ increases from $0.005$ to $0.020$, $\hat b_1$ shifts leftward by $0.1340$. This leftward shift of $\hat b_1(\sigma)$ implies that the structural transition of the transition layer occurs at lower $b_1$ values, necessitating earlier intervention to prevent HABs.

\subsection{Comparison with time series indicators}\label{subsec:classical_ews}
We compare $\mathrm{EWS}_{\mathrm{geom}}$ with two classical time series indicators, log-variance and lag-one autocorrelation. Time series statistics are derived from simulations of Eqs.~\eqref{eq:sde_T}--\eqref{eq:sde_u}. For comparison, all indicator curves are linearly normalised, detailed data preprocessing procedures are provided in~\ref{app:classical_ews}.

\begin{figure}[t]
\centering
\includegraphics[width=0.98\linewidth]{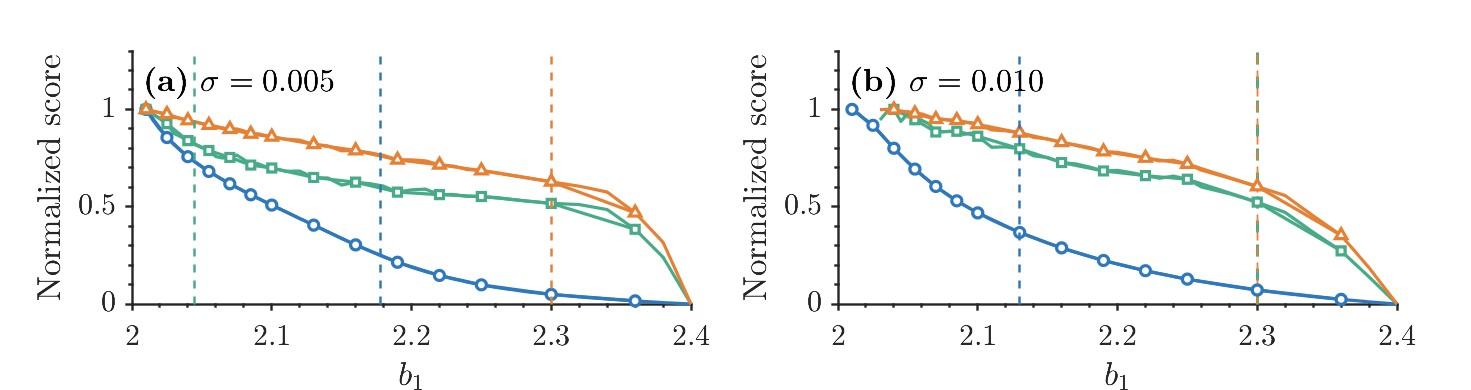}
\includegraphics[width=0.98\linewidth]{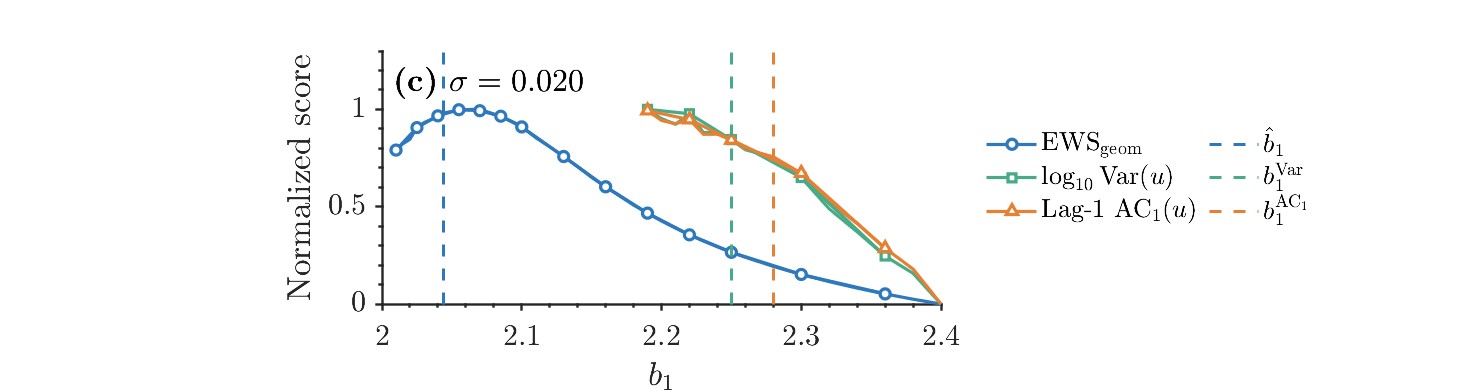}
%\caption{Comparison with time series indicators. Curves show normalised $EWS_{\mathrm{geom}}$, $\log_{10}\mathrm{Var}(u)$, and lag-one $\mathrm{AC}_1(u)$ versus $b_1$ for $\sigma = 0.005, 0.010, 0.020$. Vertical lines: BIC breakpoints ($\hat{b}_1$).}
\caption{
Comparison of the geometric metric with classical time-series indicators for (a) $\sigma = 0.005$, (b) $\sigma = 0.010$, and (c) $\sigma = 0.020$. Curves with markers represent the normalised scores of $\mathrm{EWS}_{\mathrm{geom}}$ (blue circles), $\log_{10}\mathrm{Var}(u)$ (green squares), and lag-one $\mathrm{AC}_1(u)$ (orange triangles) as functions of $b_1$. 
Vertical dashed lines indicate the corresponding critical thresholds ($\hat{b}_1$, $b_1^{Var}$, and $b_1^{AC_1}$) identified via BIC breakpoints for each indicator. 
All scores are normalised to the range $[0, 1]$ for comparison.}
\label{fig:ews_compare}
\end{figure}

Under weak noise ($\sigma = 0.005, 0.010$), $\mathrm{EWS}_{\mathrm{geom}}$ decays smoothly as $b_1$ increases, whereas time series indicators remain near baseline values until close to the bifurcation point in Fig.~\ref{fig:ews_compare}(a)--(b). This difference arises because variance and autocorrelation rely on CSD, which requires large enough fluctuations near the bifurcation to produce a signal. In contrast, $EWS_{\mathrm{geom}}$ measures the global geometric distance to the stochastic separatrix, detecting the shrinking basin volume before statistical instability becomes apparent.

Under stronger noise ($\sigma = 0.020$, Fig.~\ref{fig:ews_compare}(c)), $\log_{10}\mathrm{Var}(u)$ and $\mathrm{AC}_1(u)$ exhibit data gaps in the lower $b_1$ region. These discontinuities occur because the noise-induced MFPT becomes shorter than the observation window, precluding the accumulation of sufficient effective samples. Conversely, $EWS_{\mathrm{geom}}$ as a static geometric measure, remains computable across the full range. Quantitatively, the BIC breakpoint for $EWS_{\mathrm{geom}}$ is identified at $\hat{b}_1 = 2.044$, appearing earlier than the breakpoints for variance ($b_1^{\mathrm{Var}} \approx 2.250$) and autocorrelation ($b_1^{\mathrm{AC}_1} \approx 2.280$). This difference reflects the fact that geometric indicators respond to the reshaping of the basin boundary, while statistical indicators require the accumulation of CSD to produce a signal.

We emphasise that the geometric indicator is fundamentally model-based.
Its computation requires knowledge (or reliable reconstruction) of the
drift and diffusion structure in phase space, since it involves solving
a backward Kolmogorov boundary value problem for the committor function.
Therefore, $\mathrm{EWS}_{\mathrm{geom}}$ is most naturally applicable in
model-driven or data-assimilative contexts rather than purely observational
settings. In experimental systems where only time series are available,
classical statistical indicators may remain more practical.

%We note that the geometric indicator requires knowledge of phase space structure, which may be challenging to obtain in experimental systems. For model-based applications where the dynamics are known or can be approximated, $EWS_{\mathrm{geom}}$ provides geometric information that does not require long stationary time series.

\section{Timescales from geometry}\label{sec:timescale}
This section establishes a quantitative relation between the MFPT and the geometric early warning indicator. The relation arises from two independent asymptotic behaviours in the weak noise regime, the $1/\sigma^2$ dependence of the logarithmic transition time and the linear $\sigma$ dependence of the geometric indicator. These two behaviours are coupled through the noise intensity $\sigma$, yielding a linear relation $\log\langle\tau\rangle \sim 1/\mathrm{EWS}_{\mathrm{geom}}^2$.

\subsection{Mean first passage time}\label{subsec:mfpt_def}
Let $X_t=(T_t,u_t)$ satisfying Eqs.~\eqref{eq:sde_T}--\eqref{eq:sde_u} with infinitesimal generator $\mathcal{L}_{b_1,\sigma}$. For fixed parameters $(b_1,\sigma)$, we define the first passage time to the bloom state neighbourhood $R_{E_3}$ as:
\[
\tau_{R_{E_3}}(x;b_1,\sigma) := \inf\{t \geq 0 : X_t \in R_{E_3} \mid X_0 = x\},
\qquad x \in \Omega \setminus R_{E_3},
\]
where $\Omega \subset \mathbb{R}^2$ is a bounded computational domain containing both attractors. The MFPT function $\tau(x;b_1,\sigma) := \mathbb{E}_x[\tau_{R_{E_3}}(x;b_1,\sigma)]$ satisfies the backward Kolmogorov equation
\begin{equation}\label{eq:mfpt_bvp}
\begin{aligned}
\mathcal{L}_{b_1,\sigma}\,\tau(x;b_1,\sigma) &= -1, && x \in \Omega \setminus R_{E_3}, \\
\tau(x;b_1,\sigma) &= 0, && x \in R_{E_3}, \\
\partial_{\mathbf{n}} \tau(x;b_1,\sigma) &= 0, && x \in \partial\Omega,
\end{aligned}
\end{equation}
where $\partial_{\mathbf{n}}$ denotes the outward normal derivative. The well-posedness of Eq.~\eqref{eq:mfpt_bvp} follows from the uniform ellipticity of $\mathcal{L}_{b_1,\sigma}$ and compatibility of the boundary conditions.

To obtain a scalar transition metric independent of the initial position, we define the spatial average over the background state neighbourhood $R_{E_1}$: 
\begin{equation}\label{eq:tau_avg}
\langle \tau \rangle_{R_{E_1}}(b_1,\sigma) := \frac{1}{|R_{E_1}|} \int_{R_{E_1}} \tau(x;b_1,\sigma) \, \mathrm{d}x.
\end{equation}
This average is robust in the weak noise regime. Provided the neighbourhood size is significantly smaller than the basin of attraction, the residence distribution of trajectories within $R_{E_1}$ converges to the local stationary density $\rho_{E_1}(x) \propto \exp(-2V(x)/\sigma^2)$. Since the quasipotential $V(x)$ varies slowly within the strictly local region $R_{E_1}$, $\rho_{E_1}(x)$ is asymptotically uniform, making the spatial average equivalent to the probabilistic average up to $O(\sigma^2)$ corrections.

\subsection{Asymptotic geometric--temporal scaling in weak-noise limit}\label{subsec:theory_link}
Consider the weak noise regime where cross-basin transitions become rare events. Freidlin--Wentzell large deviation theory gives the dominant order behaviour of the MFPT~\cite{FreidlinWentzell2012}:
\begin{equation}\label{eq:fw_asymp}
\log \langle \tau \rangle_{R_{E_1}}(b_1,\sigma) = \frac{\Delta(b_1)}{\sigma^2} + o(\sigma^{-2}),
\end{equation}
where $\Delta(b_1) > 0$ represents the quasipotential barrier height separating $R_{E_1}$ from $R_{E_3}$.

In the weak-noise regime the committor function develops a sharp
transition layer in the normal direction to the stochastic separatrix
$\Gamma$. Introducing a local normal coordinate $\xi$ and retaining the
leading-order terms yields a one-dimensional boundary-layer equation
for $q$ in the normal direction. The resulting solution has the classical error–function profile
characteristic of diffusive boundary layers arising in singularly
perturbed advection–diffusion equations. 

The detailed boundary-layer analysis leading to this result is
presented in~\ref{Deriv}. The resulting
local transition-layer width between the isocommittor surfaces
$q=1/2\pm\alpha$ is
\[
w_\alpha(s;b_1,\sigma)
= C_\alpha\,\sigma\sqrt{\frac{\widetilde a_{nn}}{\lambda}},
\]

where $\lambda$ denotes the local repulsion rate of the separatrix
and $\tilde a_{nn}$ is the normal component of the diffusion tensor.
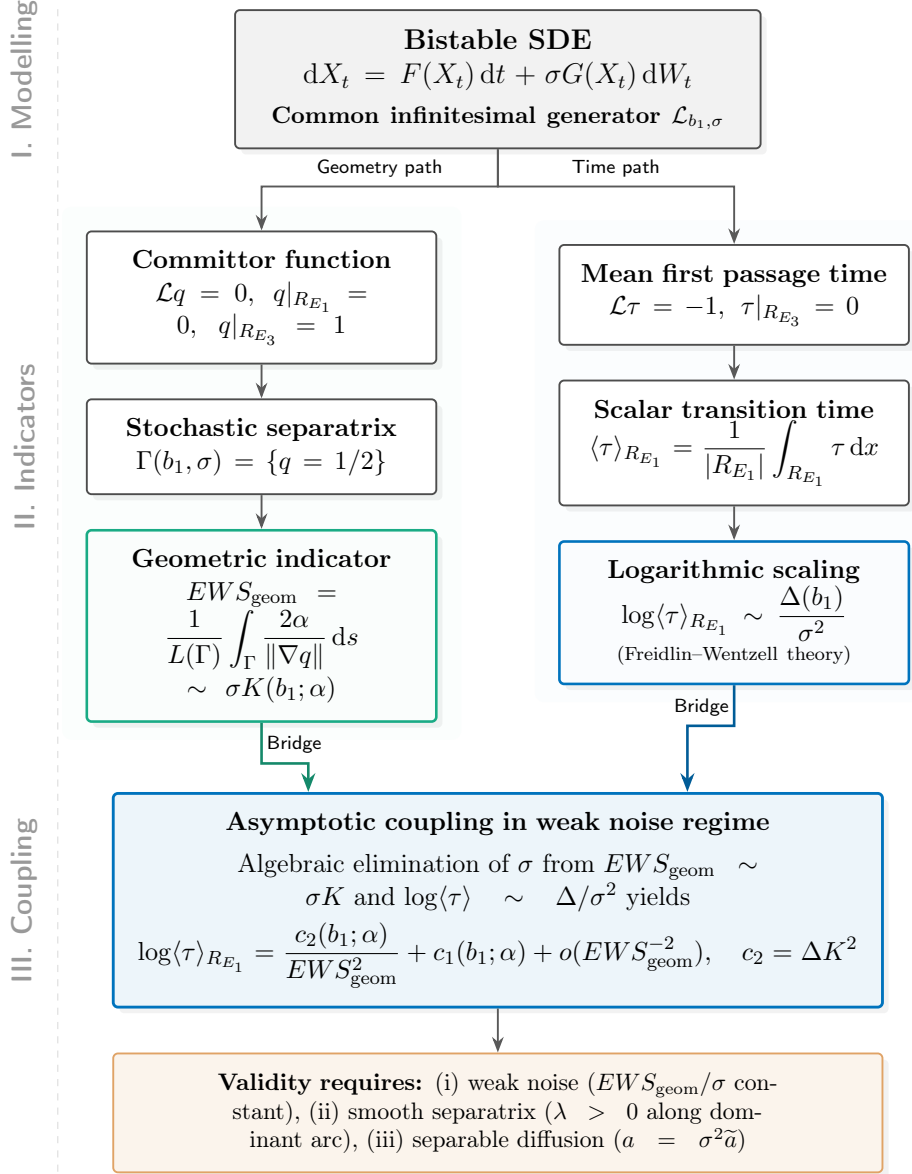
\begin{figure}[htbp]
\vspace{-2mm}
\centering
\scalebox{0.9}{
\begin{tikzpicture}[
    node distance=0.5cm and 0.9cm,
    >=Stealth,
    base/.style={
        rectangle, rounded corners=2pt, draw=black!70, thick, 
        fill=white, align=center, inner sep=8pt, text width=4.6cm,
        font=\small, drop shadow={opacity=0.1}
    },
    root/.style={base, fill=gray!10, text width=7.2cm, font=\bfseries},
    highlightGreen/.style={base, draw=myGreen, line width=1.1pt, fill=bgGreen!30},
    highlightBlue/.style={base, draw=myBlue, line width=1.1pt, fill=bgBlue!30},
    relation/.style={base, draw=myBlue, fill=bgBlue, text width=10.8cm, line width=1.2pt, font=\small},
    condition/.style={base, draw=myOrange!80, fill=bgOrange, text width=10.8cm, font=\footnotesize},
    phase/.style={font=\sffamily\bfseries\color{gray!80}, rotate=90, anchor=center}
]

    % --- Phase I: Modeling ---
    \node[root] (sde) {
        Bistable SDE \\
        $\mathrm{d}X_t = F(X_t)\,\mathrm{d}t + \sigma G(X_t)\,\mathrm{d}W_t$ \\
        \vspace{2pt}
        \footnotesize{Common infinitesimal generator $\mathcal{L}_{b_1,\sigma}$}
    };

    % --- Phase II: Geometry Path ---
    \node[base, below left=1.2cm and -3cm of sde] (q) {
        \textbf{Committor function} \\
        $\mathcal{L}q = 0,\; q|_{R_{E_1}}=0,\; q|_{R_{E_3}}=1$
    };
    \node[base, below=of q] (gamma) {
        \textbf{Stochastic separatrix} \\
        $\Gamma(b_1,\sigma) = \{q = 1/2\}$
    };
    \node[highlightGreen, below=of gamma] (ews) {
        \textbf{Geometric indicator} \\
        $\displaystyle EWS_{\mathrm{geom}} = \frac{1}{L(\Gamma)}\int_{\Gamma}\frac{2\alpha}{\|\nabla q\|}\,\mathrm{d}s$ \\
        \vspace{2pt}
        $\sim \sigma K(b_1;\alpha)$
    };

    % --- Phase II: Time Path ---
    \node[base, below right=1.4cm and -3cm of sde] (tau) {
        \textbf{Mean first passage time} \\
        $\mathcal{L}\tau = -1,\; \tau|_{R_{E_3}}=0$
    };
    \node[base, below=of tau] (mfpt) {
        \textbf{Scalar transition time} \\
        $\displaystyle \langle\tau\rangle_{R_{E_1}} = \frac{1}{|R_{E_1}|}\int_{R_{E_1}}\tau\,\mathrm{d}x$
    };
    \node[highlightBlue, below=of mfpt] (logtau) {
        \textbf{Logarithmic scaling} \\
        $\displaystyle \log\langle\tau\rangle_{R_{E_1}} \sim \frac{\Delta(b_1)}{\sigma^2}$ \\
        \vspace{2pt}
        \scriptsize{(Freidlin--Wentzell theory)}
    };

    % --- Phase III: Coupling ---
    \node[relation, below=1.3cm of $(ews.south)!0.5!(logtau.south)$] (couple) {
        \textbf{Asymptotic coupling in weak noise regime} \\[4pt]
        Algebraic elimination of $\sigma$ from $EWS_{\mathrm{geom}} \sim \sigma K$ and $\log\langle\tau\rangle \sim \Delta/\sigma^2$ yields \\[4pt]
        $\displaystyle \log\langle\tau\rangle_{R_{E_1}} = \frac{c_2(b_1;\alpha)}{EWS_{\mathrm{geom}}^2} + c_1(b_1;\alpha) + o(EWS_{\mathrm{geom}}^{-2}), \quad c_2 = \Delta K^2$
    };

    \node[condition, below=0.65cm of couple] (cond) {
        \textbf{Validity requires:} (i) weak noise ($EWS_{\mathrm{geom}}/\sigma$ constant), 
        (ii) smooth separatrix ($\lambda > 0$ along dominant arc), (iii) separable diffusion ($a = \sigma^2 \widetilde{a}$)
    };

    % --- Connectors ---
    \begin{scope}[arrows={-Stealth[scale=1.0]}, color=black!65, thick]
        \draw (sde.south) -- ++(0,-0.55) -| (q.north) node[pos=0.25, above, font=\scriptsize\sffamily, black] {Geometry path};
        \draw (sde.south) -- ++(0,-0.55) -| (tau.north) node[pos=0.25, above, font=\scriptsize\sffamily, black] {Time path};
        
        \draw (q) -- (gamma);
        \draw (gamma) -- (ews);
        \draw (tau) -- (mfpt);
        \draw (mfpt) -- (logtau);
        
        \draw[myGreen!80!black, line width=1.1pt] (ews.south) -- ++(0,-0.6) -| ([xshift=-2.8cm]couple.north) 
            node[pos=0.35, above, black, font=\scriptsize\sffamily] {Bridge};
        \draw[myBlue!80!black, line width=1.1pt] (logtau.south) -- ++(0,-0.6) -| ([xshift=2.8cm]couple.north) 
            node[pos=0.35, above, black, font=\scriptsize\sffamily] {Bridge};
        \draw (couple) -- (cond);
    \end{scope}

    % --- Background & Phase Labels ---
    \begin{scope}[on background layer]
        \coordinate (lineX) at (-6.5, 0);
        \draw[dashed, gray!30] (lineX |- sde.north) -- (lineX |- cond.south);
        
        \node[phase] at (-7.0, 0 |- sde.center) {I. Modelling};
        \node[phase] at (-7.0, 0 |- gamma.center) {II. Indicators};
        \node[phase] at (-7.0, 0 |- couple.center) {III. Coupling};

        \fill[bgGreen!25, rounded corners=4pt] ($(q.north west)+(-0.35,0.35)$) rectangle ($(ews.south east)+(0.35,-0.25)$);
        \fill[bgBlue!25, rounded corners=4pt] ($(tau.north west)+(-0.35,0.35)$) rectangle ($(logtau.south east)+(0.35,-0.25)$);
    \end{scope}

\end{tikzpicture}
}
\caption{
Schematic overview of the geometric--temporal coupling mechanism.
The upper branch (geometry path) derives the stochastic separatrix and
the geometric indicator $EWS_{\mathrm{geom}}$ from the committor equation.
The lower branch (time path) derives the mean first passage time scaling
via Freidlin--Wentzell theory.
Elimination of $\sigma$ yields the asymptotic relation
$\log\langle\tau\rangle \sim 1/EWS_{\mathrm{geom}}^2$
under weak-noise conditions.
}
\label{fig:schematic}
\end{figure}

Averaging along the arc length of $\Gamma$ gives the geometric
indicator

\begin{equation}\label{eq:ews_asymp}
EWS_{\mathrm{geom}}(b_1,\sigma;\alpha)
= \sigma K(b_1;\alpha) + o(\sigma),
\end{equation}
where $K(b_1;\alpha)>0$ is a geometric coefficient obtained by
averaging the local diffusion–repulsion balance along the
separatrix (\ref{Deriv}).

Combining the Freidlin--Wentzell scaling Eq.~\eqref{eq:fw_asymp} with the geometric
width scaling Eq.~\eqref{eq:ews_asymp} yields the main geometric--temporal coupling
relation
\begin{equation}\label{eq:logtau_ews_relation}
\log \langle \tau \rangle_{R_{E_1}}(b_1,\sigma) = \frac{c_2(b_1;\alpha)}{\mathrm{EWS}_{\mathrm{geom}}(b_1,\sigma;\alpha)^2} + c_1(b_1;\alpha) + o\bigl(\mathrm{EWS}_{\mathrm{geom}}^{-2}\bigr)
\end{equation}
which establishes an explicit asymptotic link between the geometry of
the stochastic separatrix and the time scale of noise-induced basin
transitions.

The slope satisfies $c_2(b_1;\alpha) = \Delta(b_1) \, K(b_1;\alpha)^2$, the intercept $c_1$ depends on the reference timescale (e.g., domain size) and lacks a universal form, whereas $c_2$ encodes the fundamental asymptotic coupling between geometry and 
timescales (Section~\ref{subsec:robustness}). The Eq.~\eqref{eq:logtau_ews_relation} is the algebraic combination of two independent asymptotic behaviours through the common variable $\sigma$. Its linear form $\log\langle\tau\rangle \sim 1/\mathrm{EWS}_{\mathrm{geom}}^2$ is determined by the asymptotic structure and holds for bistable diffusion models satisfying the conditions of weak noise, smooth separatrix, and separable diffusion (Section~\ref{subsec:validity}). A conceptual summary of the geometric--temporal coupling mechanism
is shown in Fig.~\ref{fig:schematic}.

We stress that Eq.~\eqref{eq:logtau_ews_relation} is not a general property of arbitrary
stochastic systems, but follows specifically from the separable diffusion
structure and non-degenerate normal diffusion along the stochastic
separatrix. Systems with state-dependent multiplicative noise that
does not factorise as $\sigma^2 \tilde a$ may exhibit modified scaling.
\subsection{Numerical verification}\label{subsec:num_check}
We fix $b_1 = 2.10$ and scan $\sigma \in [0.005, 0.025]$. For each $\sigma$, $\langle \tau \rangle_{R_{E_1}}$ is computed via finite differences method and Monte Carlo simulation, relevant settings are provided in~\ref{app:mfpt}.

\begin{figure}[t]
\centering
\includegraphics[width=0.8\linewidth]{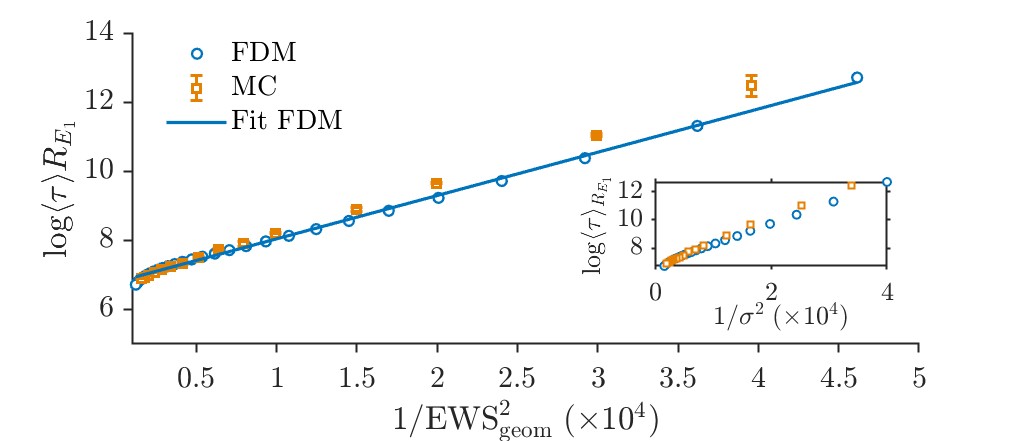}
\caption{Numerical verification of the geometric timescale relation at fixed $b_1=2.10$. Main panel: $\log\langle\tau\rangle_{R_{E_1}}$ versus $1/EWS_{geom}^2$ obtained from the finite difference method (FDM) (blue circles) and Monte Carlo (MC) simulations (orange squares), together with the linear fit to the FDM data (solid blue line). Inset: same data versus $1/\sigma^2$ ($R^2=0.9991$).}
\label{fig:mfpt_ews_fit}
\end{figure}

The Fig.~\ref{fig:mfpt_ews_fit} shows $\log \langle \tau \rangle_{R_{E_1}}$ versus $1/EWS_{\mathrm{geom}}^2$. In the weak noise regime, linear regression yields $R^2 = 0.9982$, confirming the dominant order linear structure predicted by Eq.~\eqref{eq:logtau_ews_relation}. The inset displays the same data versus $1/\sigma^2$, yielding $R^2 = 0.9991$ and validating the Freidlin--Wentzell asymptotic relation in Eq.~\eqref{eq:fw_asymp}. Monte Carlo and finite differences method results agree within statistical uncertainty, confirming numerical reliability.

\subsection{Robustness analysis}\label{subsec:robustness}
\begin{figure}[t]
\centering
\includegraphics[width=0.98\linewidth]{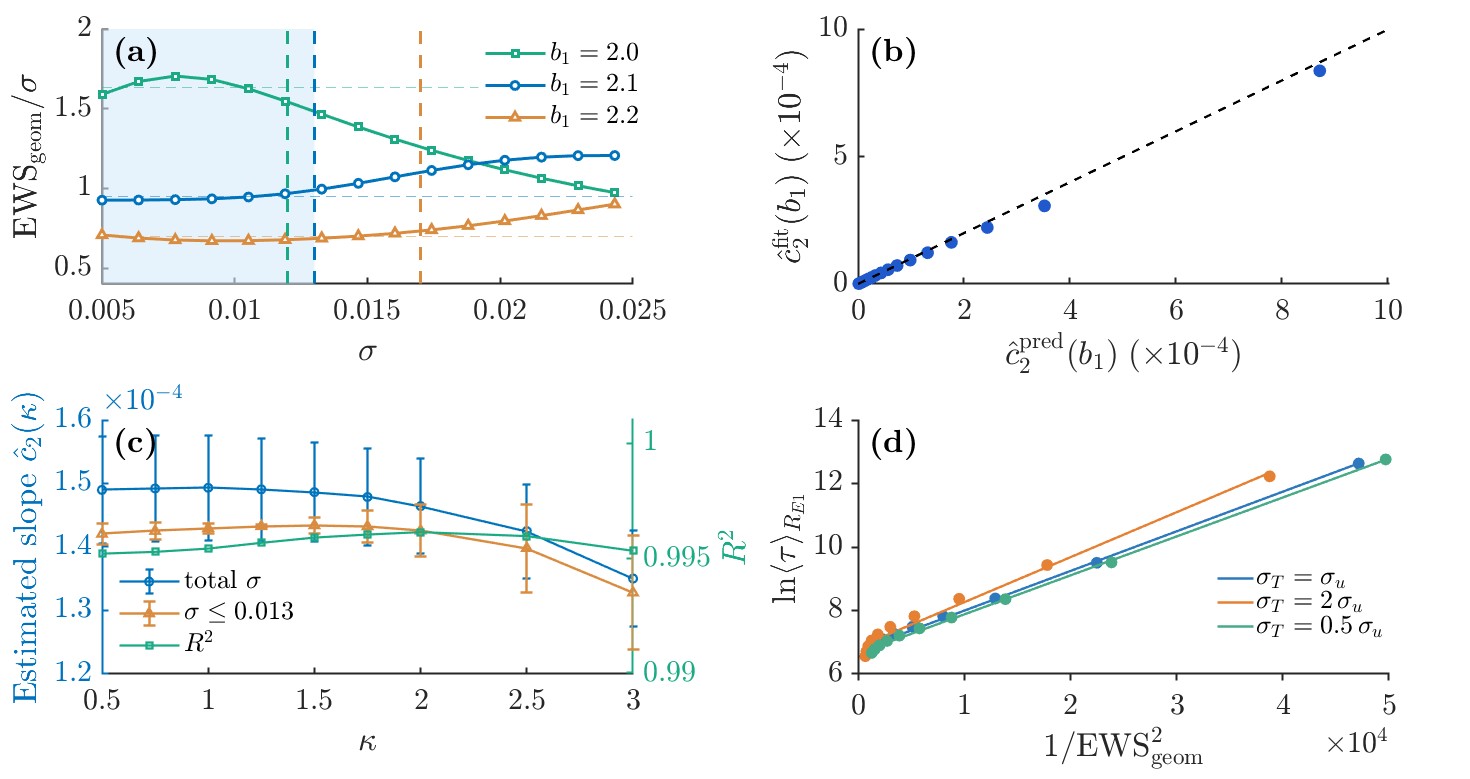}
\caption{Robustness tests for the $\log\langle\tau\rangle$--$1/EWS_{\mathrm{geom}}^2$ relation. (a) Ratio $EWS_{\mathrm{geom}}/\sigma$ versus $\sigma$ for $b_1 = 2.0$ (green), $2.1$ (blue), and $2.2$ (orange). Shaded region and vertical dashed lines denote weak noise limits for $b_1=2.10$, $\sigma=0.013$. (b) Fitted coefficient $c_2^{\mathrm{fit}}$ versus analytical prediction $c_2^{\mathrm{pred}} = \Delta K^2$ (dashed line: $y = x$, relative error $=4.6\%$). (c) Sensitivity of $\hat{c}_2$ to scaling factor where $r=\kappa r_0$ for $\kappa \in [0.5, 3.0]$, showing the estimated slope $\hat{c}_2$ evaluated using the total $\sigma$ range (blue) and the weak noise regime ($\sigma \le 0.013$) (orange), alongside the corresponding $R^2$ values (green, right axis). (d) Scaling relation for different noise ratios: $\sigma_T = \sigma_u$ (blue), $\sigma_T = 2\sigma_u$ (orange), and $\sigma_T = 0.5\sigma_u$ (green), with $R^2 \ge 0.985$}

\label{fig:timescale_robustness}
\end{figure}

The Fig.~\ref{fig:timescale_robustness} demonstrates that the affine relation in Eq.~\eqref{eq:logtau_ews_relation} arises from the separable diffusion structure $a=\sigma^2\widetilde{a}$ rather than numerical artefacts.

Panel (a) establishes the weak noise regime where $EWS_{\mathrm{geom}}$ scales linearly with $\sigma$. For $\sigma\leq0.013$ the ratio $EWS_{\mathrm{geom}}/\sigma$ maintains relative fluctuation below $5\%$, validating the asymptotic scaling in Eq.~\eqref{eq:ews_asymp}; beyond $\sigma=0.017$, the deviation from linearity increases as the one-dimensional normal approximation fails.

Panel (b) validates the theoretical prediction $c_2=\Delta K^2$. By independently estimating $\Delta$ (from $\log\langle\tau\rangle$ vs $1/\sigma^2$) and $K$ (from $EWS_{\mathrm{geom}}/\sigma$), we find that the fitted coefficients agree with the prediction typically within $5\%$ relative error. Larger deviations occur only near the bifurcation boundaries where the vanishing normal repulsion $\lambda$ challenges the asymptotic assumptions.

Panel (c) confirms insensitivity to the basin definition. Varying the neighbourhood size ($r = \kappa r_0$) over the range $\kappa \in [0.5, 3.0]$ alters the fitted coefficient by less than $5\%$ while maintaining $R^2 > 0.995$.

Panel (d) examines asymmetric noise structures ($\sigma_T=2\sigma_u$ and $\sigma_T=0.5\sigma_u$). The linear form persists ($R^2>0.985$) though coefficients vary, reflecting the geometric dependence of $K$ on the diffusion matrix $\widetilde{a}$ while the coupling form remains universal.

~\ref{app:robustness} reports additional tests on grid refinement and domain size. All variations change the coefficient by less than $4\%$ with $R^2\gtrsim0.997$, confirming that Eq.~\eqref{eq:logtau_ews_relation} is a structural consequence of separable noise intensity.

\subsection{Validity conditions}\label{subsec:validity}

The correspondence between the geometric indicator and the transition timescale depends on three asymptotic conditions derived in Section~\ref{subsec:robustness}.

First, the noise intensity should remain within the weak noise regime where $\mathrm{EWS}_{\mathrm{geom}} \propto \sigma$. We define validity by relative fluctuation of $\mathrm{EWS}_{\mathrm{geom}}/\sigma < 5\%$. As shown in Fig.~\ref{fig:timescale_robustness}(a), at fixed $b_1=2.1$, the linear regime holds up to $\sigma \approx 0.013$. 

Second, the normal repulsion rate $\lambda$ should remain positive along the stochastic separatrix. Near bifurcation points (e.g., $b_1=2.0$), the vanishing normal repulsion ($\lambda \to 0$) leads to deviations in the geometric scaling, as observed in the numerical results. When $\lambda$ becomes too small, the one-dimensional normal approximation breaks down.

Third, the diffusion tensor should maintain the separable form $a_{ij} = \sigma^2 \widetilde{a}_{ij}$. This ensures that the noise magnitude $\sigma$ factors out from the geometric integral $K$, allowing separation of variables in Eq.~\eqref{eq:logtau_ews_relation}.

To verify that these conditions apply beyond the specific ecological model, we analysed the one-dimensional Schl\"ogl model (\ref{app:schlogl}). The results confirm that the linear relation $\log\langle\tau\rangle \sim 1/\mathrm{EWS}_{\mathrm{geom}}^2$ holds with high precision within the identified validity domain.

To summarise, the admissible validity domain can be characterised
quantitatively by: (i) relative deviation $\mathrm{EWS}_{\mathrm{geom}}/\sigma < 5\%$,
(ii) strictly positive normal repulsion $\lambda > 0$ along $\Gamma$,
and (iii) separable diffusion structure $a = \sigma^2 \tilde a$.
As the system approaches criticality, the admissible noise interval
naturally contracts due to curvature growth and weakening normal
repulsion, limiting the asymptotic geometric approximation.

\section{Conclusion and future directions}\label{sec:conclusion}

This study addresses a fundamental challenge in detecting EWS for Arctic UIBs: the inability of short, irregular observational records to yield reliable statistical estimates such as variance and autocorrelation. By analysing a stochastic temperature--phytoplankton model exhibiting bistability, we identified a geometric precursor based on the probabilistic transition layer. The committor function defines a stochastic separatrix $\Gamma=\{q=1/2\}$, and the normal width of its transition layer yields a geometric indicator, $\mathrm{EWS}_{\mathrm{geom}}$, which quantifies the stability of the basin boundary using phase space geometry alone.

Our analysis reveals that the loss of basin stability involves two decoupled mechanisms: (i) the geometric shift of the separatrix named MDB and MDS; (ii) the noise-induced broadening of the transition layer $\mathrm{EWS}_{\mathrm{geom}}$. Crucially, the latter provides an earlier warning signal. The BIC breakpoint for $\mathrm{EWS}_{\mathrm{geom}}$ occurs at $\hat b_1=2.044$, significantly preceding the breakpoints for variance ($b_1^{\mathrm{Var}}=2.250$) and lag-one autocorrelation ($b_1^{\mathrm{AC_1}}=2.280$). This lead arises because $\mathrm{EWS}_{\mathrm{geom}}$ responds directly to the geometric contraction of the basin of attraction, whereas classical indicators rely on CSD, which requires a more advanced stage of destabilisation to manifest a detectable statistical signal.

From a dynamical systems perspective, the geometric indicator $\mathrm{EWS}_{\mathrm{geom}}$ can also be interpreted as a measure of system resilience. A narrow transition layer corresponds to a sharply defined basin boundary and therefore to strong resilience of the attractor against stochastic perturbations, whereas a broad transition layer indicates reduced resilience and increased susceptibility to noise-induced basin switching.

In the weak noise regime, we derived a robust scaling law linking the geometric and temporal pathways. Eliminating the noise intensity $\sigma$ from the asymptotic scalings $\mathrm{EWS}_{\mathrm{geom}}\sim\sigma$ and $\log\langle\tau\rangle_{\mathrm{RE}_1}\sim1/\sigma^2$ yields the affine relation $\log\langle\tau\rangle_{\mathrm{RE}_1}=c_1+c_2/\mathrm{EWS}_{\mathrm{geom}}^2$. This linear structure persists under variations in numerical discretisation, neighbourhood size ($R^2>0.995$), and varying noise ratios, and is successfully reproduced in the canonical Schl\"ogl model (\ref{app:schlogl}). 
These results indicate that while the coefficient $c_2=\Delta K^2$ varies across models depending on the potential barrier height $\Delta$ and local separatrix geometry $K$, 
%the linear form remains universal.
the linear form holds within the asymptotic regime characterised by
weak noise, smooth separatrix geometry, and separable diffusion structure.
Outside this domain, particularly near bifurcation points where the
normal repulsion $\lambda \to 0$, the asymptotic reduction may fail,
and deviations from linearity are expected.

The geometric indicator $EWS_{geom}$ offers distinct practical advantages in strong noise or data-limited regimes. For instance, under large noise, MFPT becomes shorter than typical observational windows, leading to trajectory escape that invalidates time-series statistics. In contrast, $\mathrm{EWS}_{\mathrm{geom}}$ is a static geometric quantity computable from a single phase space snapshot. This feature directly addresses Arctic observational constraints: even limited profile data (e.g., discrete conductivity, temperature, and depth fluorescence casts from N-ICE2015) can theoretically yield a geometric assessment of transition risk without the need for long, stationary time series.

We conclude that geometric and temporal pathways serve as complementary tools. The geometric indicator is applicable when phase space structure is accessible via models or data assimilation, whereas time series indicators remain appropriate for purely observational settings with limited state information. The validity domain of the geometric approach is defined by weak noise (relative deviation of $\mathrm{EWS}_{\mathrm{geom}}/\sigma < 5\%$), smooth separatrix geometry ($\lambda>0$), and separable diffusion structure ($a=\sigma^2\widetilde{a}$). As the system approaches criticality, the admissible noise range naturally contracts, reflecting the breakdown of the geometric approximation near bifurcation points.

Future work will focus on bridging the gap between this theoretical framework and operational constraints. First, as $\lambda\to0$ near critical parameters, the separatrix curvature increases substantially; developing curvature-adaptive width definitions could extend the validity of geometric warning closer to bifurcation points. Second, since real observations typically provide only partial state variables, combining $\mathrm{EWS}_{\mathrm{geom}}$ with delay coordinate embedding or data assimilation could allow for the reconstruction of phase space geometry from limited trajectories. 

Third, addressing the intermittent nature of light penetration in the Arctic, which generates coloured rather than white noise, will be essential to adapt the $\mathrm{EWS}_{\mathrm{geom}}\propto\sigma$ scaling to realistic environmental forcing conditions. 

Finally, an especially interesting direction concerns situations where the deterministic separatrix is not a smooth invariant manifold but instead possesses an intrinsic geometric width or a fractal-like structure, as discussed for example in~\cite{LucariniBodai2020,Mehlingetal2024}. In such systems the basin boundary may already exhibit a finite thickness even in the deterministic limit. Investigating whether $\mathrm{EWS}_{\mathrm{geom}}$ is able to detect this intrinsic boundary width, and how this interacts with noise-induced broadening of the transition layer, represents an important direction for future research.

From a biological standpoint, the geometric indicator translates
stochastic basin deformation into a measurable precursor of bloom onset.
In highly variable Arctic systems, where rapid transitions limit the
reliability of conventional EWS, such a structural
indicator may improve the anticipation of abrupt biomass shifts.
Although its computation requires an underlying dynamical model,
the approach provides a mechanistic framework for integrating
process-based modelling with ecological monitoring strategies.

\section*{Glossary of abbreviations}

For ease of reference, we collect here the main abbreviations used throughout the manuscript.

\begin{description}

\item[UIBs] Under-ice blooms: rapid phytoplankton biomass accumulations occurring beneath Arctic sea ice under favourable light and stratification conditions.

\item[HABs] Harmful algal blooms: proliferations of phytoplankton species that produce toxins or otherwise disrupt marine ecosystems and food webs.

\item[SDE] Stochastic differential equation.

\item[FPE] Fokker--Planck equation associated with the SDE.

%\item[SBD] Stochastic bifurcation diagram (expectation-based representation of stationary statistics under noise).

\item[CSD] Critical slowing down, referring to the reduction of recovery rate near deterministic bifurcation points.

\item[EWS] Early warning signal.

\item[$EWS_{\mathrm{geom}}$] Geometric early warning indicator defined as the arc-length averaged transition-layer width around the stochastic separatrix (Eq.~\eqref{eq:local_width}).

\item[MFPT] Mean first passage time to a target set (Section 5.1).

\item[$\langle \tau \rangle_{R_{E_1}}$] Spatially averaged MFPT over the background-state neighbourhood $R_{E_1}$.

\item[$R_{E_1}, R_{E_3}$] Neighbourhoods of the background equilibrium $E_1$ and bloom equilibrium $E_3$, respectively.

\item[$\Gamma$] Stochastic separatrix defined as the $1/2$-isocommittor surface $\{ q = 1/2 \}$.

\item[MDB] Mean distance to the deterministic separatrix $\Gamma_{\mathrm{det}}$ (Eq.~\eqref{eq:MDB_def}).

\item[MDS] Mean distance to a reference stochastic separatrix (Eq.~\eqref{eq:MDS_def}).

\item[BIC] Bayesian Information Criterion used for breakpoint detection.

%\item[FDM] Finite difference method.
%\item[MC] Monte Carlo simulation.
\item[$\Delta(b_1)$] Quasipotential barrier height in the Freidlin--Wentzell asymptotic regime.

\item[$\lambda$] Normal repulsion rate along the stochastic separatrix (Section 5.2).

\end{description}

\section*{CRediT authorship contribution statement}
Yuzhu Shi: Writing – original draft, Writing – review \& editing, 
Visualisation, Software, Methodology, Investigation, Formal analysis.

Larissa Serdukova: Writing – review \& editing, Supervision, 
Investigation, Formal analysis, Methodology, Conceptualisation.

Yayun Zheng: Writing – review \& editing, Investigation, Conceptualisation.

Sergei Petrovskii: Writing – review \& editing, Funding acquisition, Supervision, 
Investigation, Formal analysis, Methodology, Conceptualisation.

Valerio Lucarini: Writing – review \& editing, Funding acquisition, 
Investigation, Formal analysis, Conceptualisation.

\section*{Data availability }
The code and numerical routines used to generate the results and figures
in this study are available from the corresponding author upon reasonable request.

\section*{Acknowledgements}
V.L. acknowledges partial support from the Horizon Europe projects
\emph{Past2Future} (Grant No.~101184070) and \emph{ClimTIP} (Grant No.~100018693),
as well as from the ARIA project SCOP-PR01-P003
\emph{Advancing Tipping Point Early Warning (AdvanTip)}.
V.L. also acknowledges support from the NNSFC International Collaboration
Fund for Creative Research Teams (Grant No.~W2541005).
S.P. and V.L. were partially supported by the European Space Agency
under contract No.~4000146344/24/I-LR.
S.P. also acknowledges support from the RUDN University Strategic
Academic Leadership Program.
Y.S. acknowledges the China Scholarship Council (CSC) for fully
supporting her PhD studies.

%\section*{Funding }
%This work was funded by the Advanced Research and Invention Agency - ARIA  under Grant RP201W844843, AdvanTip: Advancing Tipping Point Early Warning.
\section*{Declaration of competing interest}
The authors declare that they have no known competing financial interests
or personal relationships that could have appeared to influence the work
reported in this paper.

\appendix
%======================================================================
\section{Numerical methods and implementation}\label{app:numerical}
%======================================================================

All computations are performed in MATLAB R2024b. Model parameters follow Table~\ref{tab:parameters}. The diffusion regularisation parameter is $\delta = 10^{-4}$.

%======================================================================
\subsection{Computational domain, grid, and elliptical neighbourhoods}\label{app:setup}
%======================================================================
Backward boundary value problems (committor and MFPT) are solved on the rectangular domain$$\Omega = [0.30, 0.60] \times [0, 0.13].$$The boundaries of $\Omega$ are defined to include all equilibrium points $(T, u)$ within the bistable regime for the parameter range $b_1 \in [1.996, 2.471]$.

Homogeneous Neumann conditions $\partial_n(\cdot) = 0$ are imposed on $\partial\Omega$ using mirror reflection: for a boundary point at $x = x_{\min}$, the ghost value satisfies $u_{-1} = u_1$, yielding first order accurate normal derivative $\partial_n u = 0$.

The main results (Figs.~\ref{fig:mfpt_ews_fit}--\ref{fig:timescale_robustness}) use a uniform grid of $141 \times 141$ points (grid spacing $h_T = 2.14 \times 10^{-3}$, $h_u = 9.29 \times 10^{-4}$). Grid convergence was verified by comparing solutions on $101 \times 101$, $141 \times 141$, and $181 \times 181$ grids; relative differences in $\log\langle\tau\rangle_{R_{E_1}}$ and $EWS_{\mathrm{geom}}$ are below $1\%$ between $141 \times 141$ and $181 \times 181$.

Equilibria $E_1(b_1)$ and $E_3(b_1)$ are computed numerically using \texttt{fsolve}. The neighbourhoods $R_{E_1}$ and $R_{E_3}$ are defined as the intersection of axis-aligned ellipses centred at the equilibria and the physical domain $\Omega$ ($u \geq 0$, i = 1,3):
\[
R_{E_i} = \left\{ (T,u) \in \Omega : \left( \frac{T - E_i^{(T)}}{r_T} \right)^2 + \left( \frac{u - E_i^{(u)}}{r_u} \right)^2 \leq 1,\quad u \geq 0 \right\}
\]
with semi-axes $r_T = 0.018$ and $r_u = 0.008$. The committor satisfies $q|_{R_{E_1}} = 0$, $q|_{R_{E_3}} = 1$; the MFPT satisfies $\tau|_{R_{E_3}} = 0$.

For sample path simulations, trajectories are confined to $\Omega$ by mirror reflection at $\partial\Omega$, followed by non negativity enforcement $u \leftarrow \max\{u, 0\}$.

%======================================================================
\subsection{Committor computation and stochastic separatrix extraction}\label{app:committor}
%======================================================================

The committor $q(T,u; b_1,\sigma)$ solves the backward Kolmogorov Eq.~\eqref{eq:committor_bvp} with generator in Eq.~\eqref{eq:generator_L}. Spatial derivatives are discretised using second order centred differences:
\[
\partial_T \varphi \approx \frac{\varphi_{i+1,j} - \varphi_{i-1,j}}{2h_T}, \quad
\partial_{TT} \varphi \approx \frac{\varphi_{i+1,j} - 2\varphi_{i,j} + \varphi_{i-1,j}}{h_T^2},
\]
with analogous formulas for $u$ derivatives. %The resulting sparse linear system is solved using MATLAB's backslash operator (\texttt{\textbackslash}).

The stochastic separatrix $\Gamma(b_1,\sigma) = \{q = 1/2\}$ is extracted using MATLAB's \texttt{contourc}function with marching squares algorithm. To reduce boundary effects on geometric quantities, contour points within relative distance $d_{\mathrm{edge}} = 0.02$ from $\partial\Omega$ are excluded:
\[
d_{\mathrm{edge}}^{(T)} = 0.02 \times (0.60 - 0.30) = 0.006, \qquad
d_{\mathrm{edge}}^{(u)} = 0.02 \times (0.13 - 0) = 0.0026.
\]
Among remaining connected components, the longest branch (measured by Euclidean arc length) is retained as $\Gamma$. Grid values of $q$ outside $[0,1]$ are truncated prior to contour extraction.

%======================================================================
\subsection{Geometric indicator $EWS_{\mathrm{geom}}$}\label{app:ews}
%======================================================================

The gradient $\nabla q$ is approximated by second order centred differences at grid centres. The gradient magnitude is regularised as
\[
\|\nabla q\| \leftarrow \max\{\|\nabla q\|, \varepsilon_\nabla\},
\qquad \varepsilon_\nabla = 10^{-6},
\]
to avoid division by zero in regions where $q$ is nearly constant. The regularised gradient is bilinearly interpolated to contour points of $\Gamma$, and the local width is computed as $w_\alpha(s) = 2\alpha / \|\nabla q(s)\|$. Arc length averaging yields $EWS_{\mathrm{geom}}$ according to Eq.~\eqref{eq:ews_geom_def}. We use $\alpha = 0.1$ throughout all results reported in the main text.

Arc-length integration employs the midpoint rule on consecutive contour points $\{x_k\}_{k=1}^M$:
\[
EWS_{\mathrm{geom}} \approx \frac{1}{L(\Gamma)} \sum_{k=1}^{M-1} w_\alpha(x_k^{\mathrm{mid}}) \, \|x_{k+1} - x_k\|,
\qquad x_k^{\mathrm{mid}} = \frac{x_k + x_{k+1}}{2}.
\]
For the $141 \times 141$ grid, \texttt{contourc} returns $M \approx 250$ points along $\Gamma$, corresponding to average spacing of $4$ grid cells. Linear interpolation to double the point density changes $EWS_{\mathrm{geom}}$ by less than $0.4\%$, confirming sufficient resolution.

%======================================================================
\subsection{Mean first passage time computation}\label{app:mfpt}
%======================================================================

The MFPT field $\tau(x; b_1,\sigma)$ solves Eq.~\eqref{eq:mfpt_bvp}. Discretisation follows the same second order centred difference scheme as the committor problem. The linear system is solved using sparse LU decomposition.

The scalar transition time $\langle\tau\rangle_{R_{E_1}}$ is obtained by arithmetic averaging over grid points whose centres lie inside $R_{E_1}$. For Monte Carlo validation, we generate $N_{\mathrm{traj}} = 10^4$ independent trajectories per $(b_1,\sigma)$ using Euler--Maruyama integration with step size $\Delta t = 10^{-2}$. Each trajectory starts from a point uniformly sampled within $R_{E_1}$ and is simulated until either:
\begin{enumerate}
\item[(i)] it enters $R_{E_3}$ (recorded as first passage time $\tau$), or
\item[(ii)] it reaches maximum simulation time $T_{\max} = 10^6$ without entering $R_{E_3}$.
\end{enumerate}
Trajectories reaching $T_{\max}$ without entering $R_{E_3}$ are handled via Kaplan--Meier tail mean estimation rather than simple discarding. For $\sigma \leq 0.0222$, the censoring fraction is below $0.3\%$, and the correction changes $\langle\tau\rangle_{R_{E_1}}$ by less than $0.1\%$. The sample mean and standard error are reported; for $\sigma = 0.020$, the relative standard error is $2.3\%$ ($95\%$ confidence level).

Time step convergence was verified by repeating simulations with $\Delta t \in \{10^{-2}, 5\times10^{-3}, 2\times10^{-3}\}$. The relative change in $\log\langle\tau\rangle_{R_{E_1}}$ is below $0.8\%$ between $\Delta t = 10^{-2}$ and $\Delta t = 2\times10^{-3}$ for $\sigma \leq 0.025$.

%======================================================================
\subsection{Classical time series indicators}\label{app:classical_ews}
%======================================================================

Time series are generated by Euler--Maruyama integration of Eqs.~\eqref{eq:sde_T}--\eqref{eq:sde_u} with step size $\Delta t = 10^{-2}$. Total simulation time is $T_{\mathrm{sim}} = 4000$ with transient discard $T_{\mathrm{tr}} = 1000$. Observations are down sampled at interval $\Delta t_{\mathrm{obs}} = 1.0$, yielding sequences $\{u_n\}_{n=1}^M$ with $M = 3000$.

For each $(b_1,\sigma)$, $N_{\mathrm{ens}} = 50$ independent trajectories are generated. Only trajectories that remain in $R_{E_1}$ throughout the observation window are retained for indicator computation (conditional sampling). For each valid trajectory, sample variance and lag-one autocorrelation are computed as described in the main text. Ensemble averages are reported.

%======================================================================
\subsection{Scaling law derivation linking geometric--temporal pathways}
\label{Deriv}
This appendix provides the detailed derivation of the geometric--temporal
scaling relation presented in Section~\ref{subsec:theory_link}. In the
weak-noise limit the committor equation develops a thin boundary layer
in the normal direction to the stochastic separatrix. The analysis below
follows a standard boundary-layer treatment of the advection--diffusion
equation governing the committor function.

In the weak-noise regime $\sigma \ll 1$, the committor function $q$
develops a sharp transition layer near the stochastic separatrix
$\Gamma=\{q=1/2\}$ whose width scales as $O(\sigma)$. Consequently,
variations of $q$ in the normal direction dominate tangential variations
along the separatrix.

Introducing a local normal coordinate $\xi$ with $\nabla q$,
the committor equation $\mathcal{L}_{b_1,\sigma} q = 0$ reduces at
leading order to the one–dimensional boundary-layer approximation

\[
b_n(\xi)\,\partial_\xi q + \frac{1}{2} a_{nn}\,\partial_{\xi\xi} q \approx 0,
\]

where $b_n(\xi)$ denotes the normal component of the drift and
$a_{nn}$ is the normal diffusion coefficient.
From Section~\ref{subsec:stoch_model}, the diffusion tensor admits
the separable structure

\[
a = \sigma^2 \widetilde a ,
\]

with $\widetilde a$ independent of $\sigma$.

Linearising the normal drift near the separatrix ($\xi=0$) yields

\[
b_n(\xi) = b_n^0 + \lambda \xi + O(\xi^2),
\qquad
b_n^0 := b_n(0), \qquad
\lambda := \partial_\xi b_n(0),
\]

where $\lambda>0$ characterises the local repulsion rate of the
separatrix, a generic geometric feature of bistable systems.
Substituting $a_{nn} = \sigma^2 \widetilde a_{nn}$ gives

\[
(b_n^0 + \lambda \xi)\,\partial_\xi q
+ \frac{1}{2}\sigma^2 \widetilde a_{nn}\,\partial_{\xi\xi} q = 0 .
\]

Multiplying by the integrating factor

\[
\exp\!\left(
\frac{2}{\sigma^2 \widetilde a_{nn}}
\left(b_n^0\xi+\tfrac12\lambda\xi^2\right)
\right)
\]

yields

\[
\partial_\xi
\left[
q
\exp\!\left(
\frac{2}{\sigma^2 \widetilde a_{nn}}
\left(b_n^0\xi+\tfrac12\lambda\xi^2\right)
\right)
\right]
=0 .
\]

Integrating once and applying the boundary conditions
$q\to0$ as $\xi\to-\infty$ and
$q\to1$ as $\xi\to+\infty$
gives the classical error–function profile

\[
q(\xi)
=
\frac12
\left[
1+
\operatorname{erf}
\!\left(
\sqrt{\frac{\lambda}{\sigma^2 \widetilde a_{nn}}}
(\xi+\xi_0)
\right)
\right],
\qquad
\xi_0 := \frac{b_n^0}{\lambda}.
\]

The parameter $\xi_0$ represents a shift in the transition centre and
does not affect the relative spacing of isocommittor surfaces.

The local transition-layer width $w_\alpha(s;b_1,\sigma)$ is defined
as the normal distance between the isosurfaces $q=1/2\pm\alpha$.
Using the monotonicity of the error function gives

\[
w_\alpha(s;b_1,\sigma)
=
\frac{2|\operatorname{erf}^{-1}(2\alpha)|}
{\sqrt{\lambda/(\sigma^2 \widetilde a_{nn})}}
=
C_\alpha \sigma
\sqrt{\frac{\widetilde a_{nn}}{\lambda}},
\qquad
C_\alpha:=2|\operatorname{erf}^{-1}(2\alpha)|.
\]

Thus the dominant-order transition-layer width is proportional to the
noise amplitude $\sigma$, with a coefficient determined by the local
balance between diffusion and repulsion near the separatrix.

Averaging along the arc length of $\Gamma$ from Eq.~\eqref{eq:ews_asymp} we have
\[
EWS_{\mathrm{geom}}(b_1,\sigma;\alpha)
=
\sigma K(b_1;\alpha)+o(\sigma).
\]
where

\[
K(b_1;\alpha)
=
\frac{C_\alpha}{L(\Gamma(b_1,\sigma))}
\int_{\Gamma(b_1,\sigma)}
\sqrt{
\frac{\widetilde a_{nn}(s,0;b_1)}
{\lambda(s;b_1)}
}
\,ds
>0 .
\]

Combining the asymptotic relations in Eq.~\eqref{eq:fw_asymp}
and Eq.~\eqref{eq:ews_asymp} yields

\[
\log \langle \tau \rangle_{R_{E_1}}
=
\frac{\Delta K^2}{EWS_{\mathrm{geom}}^2}
+
o(EWS_{\mathrm{geom}}^{-2}).
\]

This expression leads directly to the geometric–temporal
coupling relation in Eq.~\eqref{eq:logtau_ews_relation}
given in the main text.

%======================================================================
\subsection{Robustness tests}\label{app:robustness}
%======================================================================

Robustness tests examine sensitivity to numerical and modelling choices. All tests retain default settings except for the varied parameter. Results are summarised in Table~\ref{tab:robustness_summary}.

\begin{table}[H]
\centering
\caption{Robustness test summary for the $\log\langle\tau\rangle$--$1/EWS_{\mathrm{geom}}^2$ relation ($b_1 = 2.10$, $\sigma \in [0.005, 0.025]$). Relative changes are computed with respect to the default configuration ($141 \times 141$ grid, $\delta_u = 10^{-4}$, $\kappa = 1.0$).}
\label{tab:robustness_summary}
\begin{tabular}{lcc}
\toprule
Test & Parameter variation & Relative change in slope (\%) \\
\midrule
Grid resolution & $101 \times 101$ & $+1.2$ \\
 & $181 \times 181$ & $-0.5$ \\
Domain size & $+10\%$ padding & $+1.3$ \\
 & $+30\%$ padding & $+2.8$ \\
Regularisation & $\delta_u = 5 \times 10^{-5}$ & $-1.8$ \\
 & $\delta_u = 2 \times 10^{-4}$ & $+3.5$ \\
Neighbourhood size & $\kappa = 0.5$ (semi axes scaled) & $-4.65$ \\
 & $\kappa = 3.0$ & $+1.64$ \\
Noise structure & $\sigma_T = 2\sigma_u$ & $+13.4$ \\
 & $\sigma_T = 0.5\sigma_u$ & $-1.8$ \\
\bottomrule
\end{tabular}
\end{table}

All tests maintain coefficient of determination $R^2 \gtrsim 0.997$, confirming robustness of the linear structure. The slope variation under neighbourhood scaling ($\eta \in [0.5, 3.0]$) is non monotonic, with maximum deviation $4.65\%$ at $\eta = 0.5$.

%======================================================================
\section{Additional results}\label{app:additional}
%======================================================================
\subsection{Second-order correction and asymmetry of the transition layer}\label{app:second_order}
The first-order definition in Eq.~\eqref{eq:local_width} assumes linear variation of $q$ along the normal direction. When the separatrix is curved, the normal variation becomes nonlinear and the two half-widths are generally unequal. To quantify this effect, we expand $q$ to second order along $n$:
\begin{equation}
\begin{gathered}
q(x + \xi n) = \frac{1}{2} + \kappa_1(x) \xi + \frac{1}{2} \kappa_2(x) \xi^2 + O(\xi^3), \\
\kappa_1(x) := \|\nabla q(x)\|, \quad
\kappa_2(x) := n(x)^\top \big(\nabla^2 q(x)\big) n(x).
\end{gathered}
\label{eq:q_normal_expansion_2nd_app}
\end{equation}
Here $\kappa_1$ is the gradient magnitude and $\kappa_2$ is the second directional derivative of $q$ along the normal, reflecting local curvature effects.

For a given threshold $\alpha$, solving $q(x + \xi_+ n) = 1/2 + \alpha$ and $q(x - \xi_- n) = 1/2 - \alpha$ yields the asymptotic expansions
\begin{equation}
\xi_\pm(x;\alpha)
= \frac{\alpha}{\kappa_1(x)} \mp \frac{\kappa_2(x)}{2\kappa_1(x)^3} \, \alpha^2 + O(\alpha^3).
\label{eq:xi_pm_2nd_app}
\end{equation}
Their difference is
\begin{equation}
\xi_+(x;\alpha) - \xi_-(x;\alpha)
= -\frac{\kappa_2(x)}{\kappa_1(x)^3} \, \alpha^2 + O(\alpha^3),
\label{eq:width_asymmetry_app}
\end{equation}
so nonzero $\kappa_2$ induces an asymmetry: the half-width increases on the locally convex side and decreases on the locally concave side.

However, the geometric indicator $EWS_{\mathrm{geom}}$ depends on the total width $w_\alpha = \xi_+ + \xi_-$. Adding the expansions cancels the second-order terms involving $\kappa_2$, giving
\begin{equation}
w_\alpha(x) = \xi_+ + \xi_-
= \frac{2\alpha}{\|\nabla q(x)\|} + O(\alpha^3).
\label{eq:width_full_2nd_app}
\end{equation}
Hence, although curvature induces half-width asymmetry, the first-order, gradient-based definition in Eq.~\eqref{eq:ews_geom_def} captures the dominant scale of the transition layer with an $O(\alpha^3)$ error. For $\alpha = 0.1$ and the parameter ranges considered here, this first-order indicator is adequate.
%======================================================================
\subsection{Bayesian Information Criterion selected breakpoints}\label{app:bic}
%======================================================================

For fixed noise intensity $\sigma$, the scalar curve $I(b_1;\sigma)$ is fitted by a continuous two segment linear model (hinge function) with breakpoint $b_1^{(k)}$. The BIC selects the optimal breakpoint $\hat{b}_1(\sigma)$, and the warning interval $[b_1^{-}(\sigma), b_1^{+}(\sigma)]$ is defined as the set of breakpoints with $\mathrm{BIC}(k) \leq \mathrm{BIC}(k^\star) + 2$.

Table~\ref{tab:warning_interval} lists the BIC selected breakpoints and warning intervals for $EWS_{\mathrm{geom}}$ corresponding to Fig.~\ref{fig:sep_shift}(b).

\begin{table}[H]
\centering
\caption{BIC selected breakpoints $\hat{b}_1(\sigma)$ and warning intervals for $EWS_{\mathrm{geom}}$.}
\label{tab:warning_interval}
\begin{tabular}{cccc}
\toprule
$\sigma$ & $\hat{b}_1(\sigma)$ & $b_1^{-}(\sigma)$ & $b_1^{+}(\sigma)$ \\
\midrule
0.005 & 2.1780 & 2.1715 & 2.1853 \\
0.010 & 2.1300 & 2.1258 & 2.1337 \\
0.020 & 2.0440 & 2.0393 & 2.0483 \\
\bottomrule
\end{tabular}
\end{table}

%======================================================================
\subsection{Schl\"ogl model verification}\label{app:schlogl}

To verify that relation in Eq.~\eqref{eq:logtau_ews_relation} arises from the asymptotic structure of bistable diffusions rather than model specific details, we apply the identical computational pipeline to the one dimensional Schl\"ogl model, a canonical bistable system in chemical kinetics~\cite{Schlogl1972}:
\[
\mathrm{d}X_t = f(X_t)\,\mathrm{d}t + \sigma\,\mathrm{d}W_t,
\qquad
f(x) = -(x - x_1)(x - x_2)(x - x_3),
\]
with $x_1 = 0.20$, $x_2 = 0.50$, $x_3 = 0.80$. The domain $[0,1]$ is discretised with $4001$ grid points. Absorbing neighbourhoods $R_{E_1}$ and $R_{E_3}$ are defined as intervals of radius $0.015$ centred at the stable equilibria.

Noise intensities $\sigma \in \{0.005, 0.0075, \dots, 0.025\}$ are tested. Linear regressions yield:
\begin{align*}
\log\langle\tau\rangle_{R_{E_1}} &= 3.94004 + 0.00404872 / \sigma^2, && R^2 = 0.999, \\
EWS_{\mathrm{geom}} &= 1.20854\,\sigma, && R^2 = 0.999, \\
\log\langle\tau\rangle_{R_{E_1}} &= 4.35055 + 0.00565261 / EWS_{\mathrm{geom}}^2, && R^2 = 0.999.
\end{align*}
The predicted coefficient $c_2^{\mathrm{pred}} = \Delta K^2 = 0.00591339$ agrees with the fitted value $c_2^{\mathrm{fit}} = 0.00565261$ within a relative error of $4.41\%$, comparable to the $4.6\%$ average error observed in the temperature--phytoplankton model. This quantitative consistency confirms that the linear form $\log\langle\tau\rangle \sim 1/EWS_{\mathrm{geom}}^2$ is universal across bistable diffusions satisfying the three validity conditions in Section~\ref{subsec:validity}, while the coefficient value $c_2$ remains model specific.
%\vspace{-1mm}

%\section{Schematic of the geometric--temporal coupling}
%\label{app:schematic}

%\clearpage
\bibliographystyle{elsarticle-num-names}
\bibliography{name.bib} 

\end{document}